\documentclass[11pt]{article}
\usepackage{makeidx}
\usepackage{graphicx}
\usepackage{latexsym}
\usepackage{amsmath}
\usepackage{amscd}
\usepackage{amssymb}
\usepackage{a4wide}
\author{St\'ephane Launois \footnote{This research was supported by a
 Marie Curie Intra-European Fellowship within the $6^{\mbox{th}}$
 European Community Framework Programme and by Leverhulme Research Interchange Grant F/00158/X}
 \\
\\
{\small{\it Laboratoire de Math\'ematiques - UMR6056, Universit\'e de Reims}}\\
\small{{\it Moulin de la Housse - BP 1039 - 51687 REIMS Cedex 2, France}}\\
\small{e-mail : stephane.launois@univ-reims.fr}}

\title{Primitive ideals and automorphism group of $U_q^+(B_2)$.}
\date{ }

\newcommand{\para}{\addtocounter{subsubsection}{1}%
{\noindent{\it \thesubsubsection.\:}}}

\newcommand{\fin}{$\blacksquare$}
\newcommand{\preuve}{\textit{Proof.} }

\newcommand{\ov}{\overline}

\newcommand{\hc}{\mathcal{H}}
\newcommand{\spec}{\mathrm{Spec}}

\newcommand{\gk}{\mathrm{GKdim}}
\newcommand{\prim}{\mathrm{Prim}}
\newcommand{\mx}{\mathrm{Max}}
\newcommand{\haut}{\mathrm{ht}}
\newcommand{\aut}{\mathrm{Aut}}

\newcommand{\eun}{e_1}
\newcommand{\edeux}{e_2}
\newcommand{\etrois}{e_3}

\newcommand{\heun}{\widehat{e_1}}
\newcommand{\hedeux}{\widehat{e_2}}
\newcommand{\hetrois}{\widehat{e_3}}

\newcommand{\teun}{\widetilde{e_1}}
\newcommand{\tedeux}{\widetilde{e_2}}
\newcommand{\tetrois}{\widetilde{e_3}}

\input{xy}
\xyoption{all}

\begin{document}
\maketitle


\newtheorem{theo}{Theorem}[section]
\newtheorem{defi}[theo]{Definition}
\newtheorem{defis}[theo]{Definitions}
\newtheorem{rem}[theo]{Remark}
\newtheorem{rems}[theo]{Remarks}
\newtheorem{nota}[theo]{Notation}
\newtheorem{notas}[theo]{Notations}
\newtheorem{prop}[theo]{Proposition}
\newtheorem{hypo}[theo]{Hypothese}
\newtheorem{lem}[theo]{Lemma}
\newtheorem{conv}[theo]{Convention}
\newtheorem{convs}[theo]{Conventions}
\newtheorem{cor}[theo]{Corollary}
\newtheorem{obs}[theo]{Observation}
\newtheorem{obss}[theo]{Observations}
\newtheorem{rap}[theo]{Recall}
\newtheorem{raps}[theo]{Rappels}

\newtheorem{paragraphe}{}[subsection]


\begin{abstract}
Let $\mathfrak{g}$ be a complex simple Lie algebra of type $B_2$ and $q$ be a non-zero complex number which is not a root of unity. In the classical case, a theorem of Dixmier asserts that the simple factor algebras of Gelfand-Kirillov dimension 2 of the positive part  $U^+(\mathfrak{g})$ of the enveloping algebra of $\mathfrak{g}$ are isomorphic to the first Weyl algebra. In order to obtain some new quantized analogues of the first Weyl algebra, we explicitly describe the prime and primitive spectra of the positive part $U_q^+(\mathfrak{g})$ of the quantized enveloping algebra of $\mathfrak{g}$ and then we study the simple factor algebras of Gelfand-Kirillov dimension 2 of $U_q^+(\mathfrak{g})$. In particular, we show that the centers of such simple factor algebras are reduced to the ground field $\mathbb{C}$ and we compute their group of invertible elements. These computations allow us to prove that the  automorphism group of $U_q^+(\mathfrak{g})$ is isomorphic to the torus $(\mathbb{C}^*)^2$, as conjectured by Andruskiewitsch and Dumas. 
\end{abstract}
$ $

\\$ $

\begin{center} \large{\textbf{Introduction}} \end{center}


Let $\mathfrak{n}$ be a finite dimensional complex nilpotent Lie algebra. The structure of the primitive factor algebras of the enveloping algebra $U(\mathfrak{n})$ of $\mathfrak{n}$ are well-known: a theorem of Dixmier asserts that these factor algebras are actually isomorphic to Weyl algebras (see \cite{Dixmier}). In contrast, the structure of the automorphism group of $U(\mathfrak{n})$ is far from being understood: even in the case of the non-abelian three dimensional nilpotent Lie algebra, there exist wild automorphisms (see \cite{alevclassique}).

In this paper, we study the quantum case. Let $q$ be a non-zero complex number which is not a root of unity. The aim of this paper is to describe the prime and primitive spectra together with the automorphism group of the positive part $U_q^+(B_2)$ of the quantized enveloping algebra of a complex Lie algebra of type $B_2$. Recall that $U_q^+(B_2)$ 
is the $\mathbb{C}$-algebra generated by two indeterminates $e_1$ and $e_2$ subject to the quantum Serre relations:
\begin{eqnarray*}
 & & e_1^2 e_2   - (q^2 + q^{-2} ) e_1 e_2 e_1 + e_2 e_1^2 = 0 \\
 & & e_2^3 e_1 -  (q^2 + 1 + q^{-2} ) e_2^2 e_1 e_2 + (q^2 + 1 + q^{-2} ) e_2 e_1 e_2^2  - e_1 e_2^3 = 0.
\end{eqnarray*}


To describe the prime and primitive spectra of $U_q^+(B_2)$, we will use the stratification theory of Goodearl and Letzter that allows the construction of a partition of these two sets by using the action of a suitable torus on $U_q^+(B_2)$. More precisely, the torus $\hc:=(\mathbb{C}^*)^2$ acts naturally by automorphisms on $U_q^+(B_2)$ via:
$$(h_1,h_2).e_i = h_i e_i \mbox{ for all } i \in \{1,2\},$$
and the stratification theory leads to a partition (called the $\hc$-stratification) of the primitive spectrum $\prim (U_q^+(B_2))$ of $U_q^+(B_2)$ in eight "$\hc$-strata" that can be described as follows. 
Set $e_3=e_1 e_2 - q^2 e_2 e_1$, $\ov{e_3}=e_1 e_2 - q^{-2} e_2 e_1$ and recall that the center of $U_q^+(B_2)$ is the polynomial ring in two variables $\mathbb{C}[z,z']$. Then the eight $\hc$-strata are:
\begin{description}
\item[\hspace{2cm} $\bullet$] $\{ \langle z-\alpha ,z'-\beta \rangle \mid \alpha,\beta \in \mathbb{C}^* \}$,
\item[\hspace{2cm} $\bullet$] $\{ \langle z-\alpha ,z' \rangle \mid \alpha \in \mathbb{C}^* \} $,
\item[\hspace{2cm} $\bullet$] $\{ \langle z ,z'-\beta \rangle \mid \beta \in \mathbb{C}^* \}$,
\item[\hspace{2cm} $\bullet$] $\{ \langle e_3 \rangle \}$,
\item[\hspace{2cm} $\bullet$] $\{ \langle \ov{e_3} \rangle \}$,
\item[\hspace{2cm} $\bullet$] $\{ \langle e_1,e_2-\beta \rangle \mid \beta \in \mathbb{C}^* \}$,
\item[\hspace{2cm} $\bullet$] $\{ \langle e_1-\alpha , e_2 \rangle \mid \alpha \in \mathbb{C}^* \}$  
\item[\hspace{2cm} $\bullet$] and $\{ \langle e_1,e_2 \rangle \}$.
\end{description}
Observe that, except for the augmentation ideal, all maximal ideals of the center $\mathbb{C}[z,z']$ of $U_q^+(B_2)$ extend to height $2$ primitive ideals of $U_q^+(B_2)$. However they do not appear in the same $\hc$-strata; in particular, those that contain $z$ and those that contain $z'$ are not in the same $\hc$-strata. This indicates that $z$ and $z'$ must be distinguished: they do not have the same status in $U_q^+(B_2)$. This is actually a natural idea since $z$ comes "directly" from the Lie algebra of type $B_2$ contrary to $z'$ which only appears at the (quantized) enveloping algebra level. This observation will play a crucial role in the calculation of the automorphism group of $U_q^+(B_2)$. More precisely, this observation will allow us to prove that the ideal generated by $ z$ is left invariant by every automorphism of $U_q^+(B_2)$. 

Next we study the simple factors of Gelfand-Kirillov dimension $2$ of $U_q^+(B_2)$. Since, in the classical case, the simple factor algebras of Gelfand-Kirillov dimension $2$ of the enveloping algebra of a Lie algebra of type $B_2$ are isomorphic to the first Weyl algebra $A_1(\mathbb{C})$ (see \cite{Dixmier}), our aim is in fact to compare the properties 
of these factor algebras with those of the first Weyl algebra
$A_1(\mathbb{C})$. In particular, we will prove that the centers of
such simple factor algebras are reduced to $\mathbb{C}$ and we will
calculate the groups of units of these algebras. This study suggests
that we distinguish between three families of such algebras:
\\\indent  $\bullet$ Those that are obtained by factorizing
$U_q^+(B_2)$ by a height 2 maximal ideal that contains $z$. Algebras
in this first family are the so-called Weyl-Hayashi algebras (see
\cite{hayashi}); they have been studied from a ring-theoritical point
of view by Alev and Dumas (see \cite{alevdumasrigidite}), Kirkman and
Small (see \cite{kirkman}), and  Malliavin (see
\cite{malliavin}). These algebras have non-trivial invertible elements
(that is, invertible elements that do not belong to $\mathbb{C}$) and
can be presented as generalized Weyl algebras (GWA for short) over a
Laurent polynomial ring in one variable (see Section \ref{GWA} for the
definition of a GWA).
\\\indent  $\bullet$ Those that are obtained by factorizing
$U_q^+(B_2)$ by a height 2 maximal ideal that contains $z'$. Algebras
in this family also have non-trivial invertible elements and can 
be presented as generalized Weyl algebras over a
Laurent polynomial ring in one variable. However they are not
isomorphic to algebras in the first family. 
\\\indent  $\bullet$ Those that are obtained by factorizing
$U_q^+(B_2)$ by a height 2 maximal ideal that does not contain $z$
nor $z'$. These algebras do not contain non-trivial invertible
elements and cannot be presented as a GWA over a polynomial ring in one variable nor a Laurent polynomial ring in one variable. Algebras in this third family provide 
some good quantized analogues of the first Weyl algebra.
\\$ $

Finally, we calculate the automorphism group of $U_q^+(B_2)$. Note that, although the structure of the automorphism group of the augmented form $\check{U}_{q} (\mathfrak{b}^{+})$, where $\mathfrak{b}^{+}$ is the Borel subalgebra of a finite dimensional complex Lie algebra $\mathfrak{g}$, has been described in \cite{fleury} in the general case, the  automorphism group $\aut (U_q^+(\mathfrak{g}))$ of $U_q^+(\mathfrak{g})$ seems to be known only when $\mathfrak{g}$ is of type $A_2$ (see \cite{alevdumasrigidite} and \cite{caldero}). Our aim here is to describe this group in the case where $\mathfrak{g}$ is of type $B_2$. The study of this automorphism group was begun by Andruskiewitsch and Dumas (see \cite{andrusdumas}) who have obtained some partial results on $\aut (U_q^+(B_2))$ by studying natural actions of this group on the center of $U_q^+(B_2)$ (which is a polynomial ring in two variables $\mathbb{C}[z,z']$) and on the prime and primitive spectra of  $U_q^+(B_2)$. In particular, they have shown the following result that will be our starting-point in the computation of the group $\aut (U_q^+(B_2))$. Denote by $\aut_z (U_q^+(B_2))$ the sub-group of $\aut (U_q^+(B_2))$ of those automorphisms of $U_q^+(B_2)$ that fix the prime ideal generated by the central element $z$. Then we have (see \cite[Proposition 3.3]{andrusdumas}):

\begin{eqnarray}
\label{intro}
\aut_z (U_q^+(B_2)) \simeq (\mathbb{C}^*)^2.
\end{eqnarray}

Concerning $\aut \left( U_q^+(B_2) \right)$ itself, Andruskiewitsch and Dumas have conjectured that the group $\aut(U_q^+(B_2))$ is also isomorphic to the torus $(\mathbb{C}^*)^2$ (see \cite[Problem 1]{andrusdumas}). The isomorphism (\ref{intro}) opens a potential route to prove this conjecture: if we can prove that every automorphism of $U_q^+(B_2)$ fixes the ideal generated by $z$, then the isomorphism (\ref{intro}) will show that $\aut (U_q^+(B_2)) \simeq (\mathbb{C}^*)^2$. This is this route that we will follow in this article. First, using the previous study of the simple factor algebras of Gelfand-Kirillov dimension 2 of $U_q^+(B_2)$, we prove that the set of those primitive ideals that contain $z$ is left invariant by every automorphism of $U_q^+(B_2)$. Next, using the fact that $U_q^+(B_2)$ is a Jacobson ring, we conclude that the ideal generated by $z$ is also invariant under every automorphism of $U_q^+(B_2)$. 
Hence we get that $\aut (U_q^+(B_2))=\aut_z (U_q^+(B_2)) \simeq (\mathbb{C}^*)^2$, as conjectured by Andruskiewitsch and Dumas. 
As a corollary, we obtain that the action of $\aut (U_q^+(B_2))$ on $\prim (U_q^+(B_2))$ has exactly $8$ orbits that we explicitly describe.

$ $

\begin{flushleft}
\textbf{Acknowledgments.}
\end{flushleft}
We thank J. Alev, F. Dumas, T.H. Lenagan and L. Richard for very helpful conversations and comments.

\newpage 
  

\section{$U_q^+(B_2)$ and some related algebras.}
$ $

 Throughout this paper, $\mathbb{C}$ denotes the field of complex numbers and $q$ is a non-zero complex number which is not a root of unity.

\subsection{Basics on $U_q^+(B_2)$.}
$ $

 In this section, we fix the notations that will be used throughout this paper. Most of these notations are taken from \cite {andrusdumas}.

\subsubsection{The algebra $U_q^+(B_2)$.} 
$ $

We denote by $U^+$, or $U_q^+(B_2)$, the quantum enveloping algebra over $\mathbb{C}$ of the nilpotent positive part of a complex simple Lie algebra of type $B_2$. Recall (see \cite{andrusdumas}, 3.1.1) that $U^+$ is the $\mathbb{C}$-algebra generated by two indeterminates $e_1$ and $e_2$ subject to the quantum Serre relations:
\begin{eqnarray}
 & & e_1^2 e_2   - (q^2 + q^{-2} ) e_1 e_2 e_1 + e_2 e_1^2 = 0 \\
 & & e_2^3 e_1 -  (q^2 + 1 + q^{-2} ) e_2^2 e_1 e_2 + (q^2 + 1 + q^{-2} ) e_2 e_1 e_2^2  - e_1 e_2^3 = 0
\end{eqnarray}

\subsubsection{$U^+$ as an iterated Ore extension.}
 \label{extOre}  
$ $

Throughout this paper, we will adopt the following notations that agree with those of \cite{andrusdumas}. We set 
\begin{eqnarray*} 
e_3 & = & e_1 e_2 - q^2 e_2 e_1 \\
 z  & =  & e_2 e_3 - q^2 e_3 e_2 \\
 \ov{e_3} & = & e_1 e_2 - q^{-2} e_2 e_1
\end{eqnarray*}
and 
\begin{eqnarray*}
z' & = & (1-q^{-4})(1-q^{-2}) e_3 e_1 e_2 + q^{-4} (1-q^{-2}) e_3^2 + 
(1-q^{-4}) ze_1 \\
 & = & (1-q^{-2}) \left( e_3 \ov{e_3} + (1+q^{-2})ze_1 \right).
\end{eqnarray*}
$ $

Recall (see \cite{andrusdumas}, 3.1.1) that the monomials $(z^i e_3^j e_1^k e_2^l)_{(i,j,k,l) \in \mathbb{N}^4}$ form a PBW-basis of $U^+$. Hence $U^+$ is the $\mathbb{C}$-algebra generated by $\eun,\edeux,\etrois, z$ with the relations :
$$\begin{array}{lll}
e_3 z = z e_3, & & \\
e_1 z = z e_1, & e_1 e_3 = q^{-2} e_3 e_1,& \\
e_2 z = z e_2, & e_2 e_3 = q^{2} e_3 e_2 + z, & e_2 e_1 = q^{-2} e_1 e_2 - q^{-2} e_3.
\end{array}$$
In other words, $U^+$ is an iterated Ore extension that we can write as follows:
$$U^+ = \mathbb{C} [z,e_3][e_1;\sigma] [e_2;\tau,\delta],$$
where $\sigma$ denotes the automorphism of $\mathbb{C}[z,e_3]$ defined by $\sigma(z)=z$ and $\sigma(e_3)=q^{-2} e_3$, 
where $\tau$ denotes the automorphism of $\mathbb{C}[z,e_3][e_1;\sigma]$ defined by $\tau(z)=z$, $\tau(e_3)=q^{2} e_3$ and $\tau (e_1) = q^{-2} e_1$, and where $\delta$ denotes the (left) $\tau$-derivation of $\mathbb{C}[z,e_3][e_1;\sigma]$ defined by $\delta(z)=0$, $\delta(e_3)=z$ and $\delta (e_1) = -q^{-2} e_3$. In particular, the algebra $U^+$ is a Noetherian domain. To be complete, let us mention that the Gelfand-Kirillov dimension $\gk (U^+)$ of $U^+$ is 4, that this algebra is catenary and that Tauvel's height formula holds in $U^+$ (see \cite[Theorem 4.8]{goodlenCAT}).

\subsubsection{The center of  $U^+$.}  
$ $

The center $Z(U^+) $ of $U^+$ is a polynomial ring in two variables. More precisely, we have (see \cite[Lemma 3.1]{andrusdumas}):
$$Z(U^+)=\mathbb{C}[z,z'].$$

\subsubsection{Some commutation relations in $U^+$.} 
$ $

The following commutation relations can be easily obtained by induction.
\\$ $
\begin{lem}
\label{relationsdansU}
For all $k \in \mathbb{N}^*$, we have:
\begin{enumerate}
\item $\displaystyle{\etrois \eun^k  =  q^{2k} \eun^k \etrois}$ and $\displaystyle{ \eun \etrois^k  =  q^{-2k} \etrois^k \eun}$.
\item $\displaystyle{\edeux \etrois^k  =  q^{2k} \etrois^k \edeux +  \frac{q^{2k}-1}{q^2-1}\etrois^{k-1}z}$.
\item $\displaystyle{\edeux \eun^k  =  q^{-2k} \eun^k \edeux - q^{-2} \frac{q^{-4k}-1}{q^{-4}-1} \etrois \eun^{k-1}}$.
\end{enumerate}
\end{lem}

\subsection{The factor algebras $U^+ / \langle z-\alpha \rangle $ ($\alpha \in \mathbb{C}$).}
$ $

This paragraph is devoted to the factor algebras $U^+ / \langle z-\alpha \rangle $ of $U^+$ ($\alpha \in \mathbb{C}$): we give a PBW basis of these algebras. Further, in the case where $\alpha =0$, the algebra obtained turns out to be isomorphic to the quantum Heisenberg algebra and so we recall some basic properties of the quantum Heisenberg algebra.

\subsubsection{$U^+ / \langle z-\alpha \rangle $ ($\alpha \in \mathbb{C}$) as an iterated Ore extension.} 
\label{quotientzalpha}
$ $

Fix $\alpha \in \mathbb{C}$. We set $B_{\alpha}:= U^+ / \langle z- \alpha \rangle$. Further, if $x \in U^+$, then $\widehat{x}$ denotes the canonical image of $x$ in $B_{\alpha}$. Since the monomials $(z^i e_3^j e_1^k e_2^l)_{(i,j,k,l) \in \mathbb{N}^4}$ form a PBW-basis of $U^+$ (see \ref{extOre}), it is easy to show that the monomials $(\hetrois^j \heun^k \hedeux^l)_{(j,k,l) \in \mathbb{N}^3}$ form a PBW-basis of $B_{\alpha}$, so that $B_{\alpha}$ is the $\mathbb{C}$-algebra generated by $\heun,\hedeux,\hetrois$ with the relations: 
$$\begin{array}{ll}
\heun \hetrois = q^{-2} \hetrois \heun,& \\
 \hedeux \hetrois = q^{2} \hetrois \hedeux  + \alpha, & \hedeux \heun = q^{-2} \heun \hedeux - q^{-2} \hetrois.
\end{array}$$
In other words, $B_{\alpha}$ is an iterated Ore extension that we can write as follows:
$$B_{\alpha} = \mathbb{C} [\hetrois][\heun;\sigma'] [\hedeux;\tau',\delta'],$$
where $\sigma'$ denotes the automorphism of $\mathbb{C}[\hetrois]$ defined by $\sigma'(\hetrois)=q^{-2} \hetrois$, 
where $\tau'$ denotes the automorphism of $\mathbb{C}[\hetrois][\heun;\sigma']$ defined by $\tau'(\hetrois)=q^{2} \hetrois$ and $\tau' (\heun) = q^{-2} \heun$, and where $\delta'$ denotes the (left) $\tau'$-derivation of $\mathbb{C}[\hetrois][\heun;\sigma']$ defined by $\delta'(\hetrois)=\alpha$ and $\delta' (\heun) = -q^{-2} \hetrois$. 
In particular, $B_{\alpha}$ is a Noetherian domain and so the ideal of $U^+$ generated by $z-\alpha$ is completely prime. 
\\$ $

In the case where $\alpha =0$, the algebra $B:=B_0=U^+ / \langle z \rangle $ is well-known since this algebra is actually isomorphic to the quantum Heisenberg algebra. In the next section, we recall some basic properties of this algebra.

\subsubsection{The quantum Heisenberg algebra.}
\label{paraheisenberg}
$ $

The quantum Heisenberg algebra, denoted by $\mathbb{H}$ or $U_q^+(A_2)$, is the quantum enveloping algebra over $\mathbb{C}$ of the nilpotent positive part of a complex simple Lie algebra of type $A_2$. Recall (see, for instance, \cite{andrusdumas}, 2.2.1) that $\mathbb{H}$ is the $\mathbb{C}$-algebra generated by two indeterminates $E_1$ and $E_2$ subject to the following relations:
\begin{eqnarray}
 & & E_1^2 E_2   - (q^2 + q^{-2} ) E_1 E_2 E_1 + E_2 E_1^2 = 0 \\
 & & E_2^2 E_1   - (q^2 + q^{-2} ) E_2 E_1 E_2 + E_1 E_2^2 = 0
\end{eqnarray}
As in \cite[2.2.1]{andrusdumas}, we set $E_3:=E_1E_2 -q^2 E_2 E_1$. Hence we have (see, for instance, \cite{alevdumasrigidite}):
\begin{itemize}
\item $\mathbb{H}$ is the iterated Ore extension over $\mathbb{C}$ generated by the indeterminates $E_1,E_3,E_2$ subject to the following relations:
$$E_3 E_1 = q^{2} E_1 E_3, \hspace{5mm} E_2 E_3 = q^{2} E_3 E_2, \hspace{5mm} E_2 E_1 = q^{-2} E_1 E_2 - q^{-2} E_3.$$
In particular, $\mathbb{H}$ is a Noetherian domain.
\item The center $Z(\mathbb{H})$ of $\mathbb{H}$ is a polynomial ring in one variable $Z(\mathbb{H})=\mathbb{C}[\Omega]$, where $\Omega$ denotes the quantum Casimir, that is:
$$\Omega = (1-q^{-4})E_3 E_1 E_2 + q^{-4} E_3^2.$$
\end{itemize}
$ $

It is well-known (see, for instance, \cite[3.2.1]{andrusdumas}) that the map $f : U^+ / \langle z \rangle \rightarrow \mathbb{H}$ defined by $f(\widehat{e_i})=E_i$ ($i \in \{1,2,3\}$) is an isomorphism of $\mathbb{C}$-algebras. Note further that $f(\widehat{z'})=(1-q^{-2}) \Omega$.


\section{Prime and primitive spectra of $U^+$.}
$ $

The aim of this paragraph is to describe explicitly the prime, primitive and maximal spectra of $U^+$. In order to do this, we use the $H$-stratification theory of Goodearl and Letzter that provides a stratification of the prime and primitive spectra of $U^+$ by considering the action of a suitable torus on this algebra.

If $I$ is a non-empty subset of an algebra $A$, we denote by $\langle I \rangle_A $ the two-sided ideal of $A$ generated by $I$. To simplify the notation, if $I$ is a non-empty subset of $U^+$, we will drop the subscript and denote by $\langle I \rangle $ the two-sided ideal of $U^+$ generated by $I$. 

If $J$ is a prime ideal in an algebra $A$, we denote by $\haut(J)$ its height.

\subsection{Prime ideals of $U^+$.}
$ $

Since $q$ is not a root of unity, it is well-known (see
\cite[Section 5]{ringel}) 
that the prime ideals of $U^+$ are completely prime. As usual, we denote by $\spec(U^+)$ the set of all (completely) prime ideals of $U^+$.

\subsubsection{$\hc$-stratification of $\spec(U^+)$.}
$ $

In order to obtain a partition of the set $\spec(U^+)$, we need to consider the following action of the torus $\hc:=(\mathbb{C}^*)^{2}$ on $U^+$. The torus $\hc$ acts naturally on $U^+$ by automorphisms via:
$$(h_1,h_2).e_i=h_i e_i \quad \mbox{ for all }i \in \{ 1,2 \}.$$
(It is easy to check that the quantum Serre relations are preserved by the group $\hc$.) Recall (see \cite[3.4.1]{andrusdumas}) that this action is rational. (We send back to \cite[II.2.]{bg} for the defintion of a rational action.) A non-zero element $x$ of $U^+$ is an {\it $\hc$-eigenvector of $U^+$} if 
$h.x \in \mathbb{C}^{*}x$ for all $h \in \hc$. An ideal $I$ of $U^+$ is {\it $\hc$-invariant} if $h.I=I$ for all $h \in \hc$. We denote by {\it $\hc$-$\spec(U^+)$} the set of all $\hc$-invariant prime ideals of $U^+$. This is a finite subset of $\spec(U^+)$ since Andruskiewitsch and Dumas have shown (see \cite[3.4.2]{andrusdumas}) :

\begin{prop}
\label{hpremier}
 $U^+$ has exactly 8 $\hc$-invariant prime ideals:  
$\langle 0 \rangle$,  $\langle z \rangle$, $\langle z' \rangle$, $\langle e_3 \rangle = \langle e_3,z,z' \rangle$, $\langle \ov{e_3} \rangle = \langle \ov{e_3},z,z' \rangle$, $\langle e_1 \rangle = \langle e_1,e_3,\ov{e_3},z,z' \rangle$, $\langle e_2 \rangle = \langle e_2,e_3,\ov{e_3},z,z' \rangle$ and $\langle e_1,e_2 \rangle = \langle e_1,e_2,e_3,\ov{e_3},z,z' \rangle$. 
\\Moreover, $\haut(\langle 0 \rangle)=0$, $\haut(\langle z \rangle)=\haut(\langle z' \rangle)=1$, $\haut(\langle e_3 \rangle)=\haut(\langle \ov{e_3} \rangle)=2$, $\haut(\langle e_1 \rangle)=\haut(\langle e_2 \rangle)=3$ and  $\haut(\langle e_1,e_2 \rangle)=4$.
\end{prop}

The action of $\hc$ on $U^+$ allows via the $\hc$-stratification
theory of Goodearl and Letzter (see \cite[II.2]{bg}) the construction
of a partition of $\spec(U^+)$ as follows. If $J$ is one of the eight
$\hc$-invariant prime ideals of $U^+$, we denote by $\spec_J(U^+)$ the
{\it $\hc$-stratum of $\spec(U^+)$ associated to $J$}. Recall that
$\spec_J(U^+):=\{ P \in \spec(U^+) \mid \bigcap_{h \in \hc} h.P=J
\}$. 
Then the $\hc$-strata $\spec_J(U^+)$ ($J \in \hc$-$\spec(U^+)$) form a partition of $\spec(U^+)$ (see \cite[3.4.1]{andrusdumas}):
$$\spec(U^+) = \bigsqcup_{J \in \hc \mbox{-}\spec(U^+)} \spec_J(U^+).$$
In order to describe the prime spectrum of $U^+$, we will now describe the $8$ $\hc$-strata of $\spec(U^+)$.

\subsubsection{Description of the $\hc$-strata of $\spec(U^+)$.}
$ $

 Among the $\hc$-invariant prime ideals of $U^+$, it is useful to distinguish those that contain $z$. There are $6$ of them: $\langle z \rangle$, $\langle e_3 \rangle$, $\langle \ov{e_3} \rangle$, $\langle e_1 \rangle$, $\langle e_2 \rangle$ and $\langle e_1,e_2 \rangle$. The $\hc$-strata corresponding to these $6$ $\hc$-invariant prime ideals are easy to compute since we dispose of an explicit description of the prime spectrum of $U^+/ \langle z \rangle \simeq \mathbb{H}$ (see \cite{malliavin}). More precisely, we deduce from this description (see \cite[Th\'eor\`eme 2.4]{malliavin}) via the isomorphism $f: U^+/ \langle z \rangle \rightarrow \mathbb{H}$ introduced in Section \ref{paraheisenberg} that the sub-poset of $\spec(U^+)$ of those primes that contains $z$ is the following.  

$$
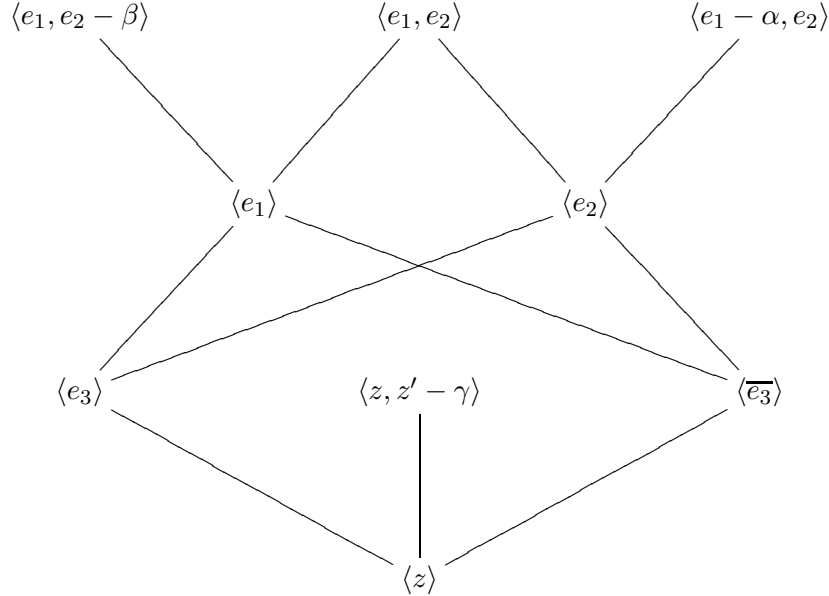
\xymatrix{ \langle e_1,e_2-\beta \rangle \ar@{-}[rdd]&& \langle  e_1,e_2 \rangle \ar@{-}[ldd]\ar@{-}[rdd] & & \langle  e_1 -\alpha , e_2 \rangle \ar@{-}[ldd] \\ 
& & & &  \\
 & \langle  e_1 \rangle \ar@{-}[ldd]\ar@{-}[rrrdd]  & &   \langle  e_2 \rangle \ar@{-}[llldd]\ar@{-}[rdd]  &  \\
& & & &  \\
\langle e_3 \rangle \ar@{-}[rrdd] & & \langle z,z'- \gamma  \rangle \ar@{-}[dd] & &  \langle  \overline{e_3} \rangle \ar@{-}[lldd]   \\
& & & &  \\
& & \langle  z \rangle   & & }$$
where $\alpha,\beta,\gamma \in \mathbb{C}^*$.\\

As a corollary, we obtain the following description of those $\hc$-strata of $\spec(U^+)$ that are associated 
to $\hc$-invariant prime ideals containing $z$.

\begin{prop}
\label{stratez}
$$\begin{array}{lll}
\spec_{\langle z \rangle } (U^+) & = & \{ \langle z \rangle  \} \cup \{ \langle z,z'-\beta \rangle \mid \beta \in \mathbb{C}^* \} ,\\
\spec_{\langle e_3 \rangle } (U^+)  & = & \{ \langle e_3 \rangle \} ,\\
\spec_{\langle \ov{e_3} \rangle } (U^+) & = & \{ \langle \ov{e_3} \rangle \},\\
\spec_{\langle e_1 \rangle } (U^+) &  = & \{ \langle e_1 \rangle  \} \cup \{ \langle e_1,e_2-\beta \rangle \mid \beta \in \mathbb{C}^* \} ,\\
\spec_{\langle e_2 \rangle } (U^+) & = & \{ \langle e_2 \rangle  \} \cup \{ \langle e_1-\alpha , e_2 \rangle \mid \alpha \in \mathbb{C}^* \} ,\\
\spec_{\langle e_1,e_2 \rangle } (U^+) & = & \{ \langle e_1,e_2 \rangle \}.
\end{array}$$
\end{prop}

It remains now to describe the $\hc$-strata associated to $\langle 0 \rangle$ and $\langle z' \rangle$. 
Let us start with the $\hc$-stratum associated to $\langle z' \rangle$. 
\\$ $

\begin{prop}
\label{stratezprime}
$\spec_{\langle z' \rangle } (U^+) = \{ \langle z' \rangle  \} \cup \{ \langle z-\alpha ,z' \rangle \mid \alpha \in \mathbb{C}^* \} $.
\\Further, for all $\alpha \in \mathbb{C}^*$, we have $\haut(\langle z-\alpha ,z' \rangle ) =2$.
\end{prop}
\preuve Observe first that the prime ideals in $\spec_{\langle z' \rangle } (U^+)$ do not contain $z$. Indeed, assume that this is not the case, that is, assume that there exists $P \in \spec_{\langle z' \rangle } (U^+)$ with $z \in P$. Then, since $z$ is an $\hc$-eigenvector, we have $z \in \bigcap_{h\in \hc} h.P = \langle z' \rangle$. This is a contradiction and so we have just proved that $\spec_{\langle z' \rangle } (U^+) \subseteq \{ P \in \spec(U^+) \mid z' \in P \mbox{ and } z \notin P \}$. On the other hand, if $P$ is a prime ideal of $U^+$ such that $z' \in P$ and  $z \notin P$, then $\bigcap_{h\in \hc} h.P $ is an $\hc$-invariant prime ideal of $U^+$ that contains $z'$ (since $z'$ is an $\hc$-eigenvector), but does not contain $z$. In view of the list of $\hc$-invariant prime ideals of $U^+$ (see Proposition \ref{hpremier}), the only possibility is $\bigcap_{h\in \hc} h.P =\langle z' \rangle $, so that 
$P \in \spec_{\langle z' \rangle } (U^+)$. To resume, we have shown that  $\spec_{\langle z' \rangle } (U^+) = \{ P \in \spec(U^+) \mid z' \in P \mbox{ and } z \notin P \}$.

We denote by $\pi$ the canonical surjection from $U^+$ onto $U^+/\langle z' \rangle$. Recall that $U^+/\langle z' \rangle$ is a Noetherian domain (see \cite[3.3.1]{andrusdumas}). Further, since $z \in Z(U^+)$, $\pi(z)$ is central in this algebra. So $\{\pi(z)^i \mid i \in \mathbb{N} \}$ is a right denominator set in $U^+/\langle z' \rangle$; the corresponding localisation of $U^+/\langle z' \rangle $ will be denoted $\displaystyle{A:=\frac{U^+}{\langle z' \rangle} [\pi(z)^{-1}] }$. Note that, since $z$ and $z'$ are $\hc$-eigenvectors, the torus $\hc$ still acts rationally by automorphisms on $A$ (see \cite[Exercise II.3.A]{bg}). Moreover, it follows from the previous study (and from classical results of non-commutative localisation theory) that the map $\varphi : P \rightarrow \displaystyle{\frac{P}{\langle z' \rangle} [\pi(z)^{-1}] }$ is an increasing bijection from $\spec_{\langle z' \rangle } (U^+)$ onto $\spec \left(A \right)$.

Before describing the prime spectrum of $A$, we establish that the algebra $A$ is $\hc$-simple in the sense of \cite[II.1.8]{bg}, that is, we show that $A$ has only one (two-sided) proper $\hc$-invariant ideal: $\langle 0 \rangle_A$. First, since $z$ and $z'$ are $\hc$-eigenvectors, we deduce from \cite[Exercise II.1.J]{bg} that the torus $\hc$ still acts by automorphisms on $A$ and that the bijection $\varphi$ induces a bijection between the set of those $\hc$-invariant prime ideals of $U^+$ that contain $z'$ but not $z$ and the set of $\hc$-invariant prime ideals of $A$. Since the set of those $\hc$-invariant prime ideals of $U^+$ that contain $z'$ but not $z$ is reduced to $\{ \langle z' \rangle \}$ (see Proposition \ref{hpremier}), we obtain that $A$ has only one $\hc$-invariant prime ideal: $\langle 0 \rangle_A$. Now, every $\hc$-invariant proper ideal of $A$ is contained in an $\hc$-invariant prime ideal of $A$, so that $\langle 0 \rangle_A$ is the unique $\hc$-invariant proper ideal of $A$. In other words, $A$ is $\hc$-simple, as desired.

We are now able to describe the prime spectrum of $A$ and the $\hc$-stratum of $U^+$ associated to $\langle z' \rangle$. Since the action of $\hc$ on $A$ is rational (see \cite[Exercise II.3.A]{bg}) and since $A$ is $\hc$-simple, it follows from  \cite[Corollary II.3.9]{bg} that extension and contraction provide mutually inverse bijections between $\spec(A)$ and $\spec(Z(A))$. Now,  in the proof of \cite[Proposition 3.4]{andrusdumas}, Andruskiewitsch and Dumas have shown that the center of $U^+/\langle z' \rangle$  is the polynomial algebra $Z(U^+/\langle z' \rangle)=\mathbb{C}[\pi(z)]$. This leads immediately to $Z(A)=\mathbb{C}[\pi(z)^{\pm 1}]$. Thus, since $\mathbb{C}$ is algebraically closed, we obtain that $\spec (Z(A))= \{ \langle 0 \rangle_{Z(A)}  \} \cup \{ \langle \pi(z)-\alpha \rangle_{Z(A)} \mid \alpha \in \mathbb{C}^* \} $ and so we deduce from the previous study that $\spec(A) = 
\{ \langle 0 \rangle_A  \} \cup \{ \langle \pi(z)-\alpha  \rangle_A \mid \alpha \in \mathbb{C}^* \} $, and next that 
$\spec_{\langle z' \rangle } (U^+) = 
\{ \varphi^{-1} (\langle 0 \rangle_A)  \} \cup \{ \varphi^{-1} (\langle \pi(z)-\alpha  \rangle_A) \mid \alpha \in \mathbb{C}^* \} $. Naturally, $\varphi^{-1} (\langle 0 \rangle_A)=\langle z' \rangle$. Thus, to achieve the proof of the first part of the proposition, it just remains to show that $\varphi^{-1} (\langle \pi(z)-\alpha  \rangle_A)=\langle z-\alpha ,z' \rangle$ for all $\alpha \in \mathbb{C}^*$. This is what we do now.

Fix $\alpha \in \mathbb{C}^*$. Clearly, we have $\varphi^{-1} (\langle \pi(z)-\alpha  \rangle_A) \supseteq \langle z-\alpha ,z' \rangle$ and so we just need to establish the reverse inclusion. Let $x \in \varphi^{-1} (\langle \pi(z)-\alpha  \rangle_A)$. Then $\pi(x)$ belongs to the ideal of $\displaystyle{A=\frac{U^+}{\langle z' \rangle} [\pi(z)^{-1}]}$  generated by $\pi(z) -\alpha$. Hence, there exist $t \in \mathbb{N}$ and $a \in U^+/\langle z' \rangle$ such that $\pi(x) \pi(z)^t=a(\pi(z) - \alpha)$. We choose such $t$ minimal. If $t> 0$, then 
$a = \frac{1}{\alpha} (a -\pi(x) \pi(z)^{t-1})\pi(z)$. Set $b:= \frac{1}{\alpha} (a -\pi(x) \pi(z)^{t-1})$. Then 
$b \in U^+/\langle z' \rangle$ and $a=b \pi(z)$. Thus, 
$\pi(x) \pi(z)^t=a(\pi(z) - \alpha)=b \pi(z) (\pi(z) - \alpha)$. Since $\pi(z) \neq 0$ and $U^+/\langle z' \rangle$ is a domain, we get that $\pi(x) \pi(z)^{t-1}=b(\pi(z) - \alpha)$ with $b \in U^+/\langle z' \rangle$. This contradicts the minimality of $t$. Hence $t=0$ and so $\pi(x) =a(\pi(z) - \alpha)$ belongs to the ideal of 
$U^+/\langle z' \rangle$ generated by $ \pi(z) -\alpha$. Thus, $x$ belongs to the ideal of $U^+$ generated by $z'$ and $z-\alpha$, as desired. This achieves the proof of the first part of the proposition.

Fix $\alpha \in \mathbb{C}^*$. It remains to prove that $\haut(\langle z-\alpha ,z' \rangle ) =2$. 
First, since $\langle 0 \rangle \varsubsetneq \langle z' \rangle \varsubsetneq \langle z-\alpha, z' \rangle $, we have  $\haut(\langle z-\alpha ,z' \rangle) \geq 2$. If $\haut(\langle z-\alpha ,z' \rangle) > 2$, then, because of the catenarity of $U^+$ (see \cite[Theorem 4.8]{goodlenCAT}), the previous chain is not saturated, so that there exists a prime ideal $P$ of $U^+$ such that either $\langle 0 \rangle \varsubsetneq P \varsubsetneq \langle z' \rangle \varsubsetneq \langle z-\alpha, z' \rangle $ or $\langle 0 \rangle \varsubsetneq \langle z' \rangle \varsubsetneq P \varsubsetneq \langle z-\alpha, z' \rangle $. Since $\haut(\langle z' \rangle) =1$, the first case cannot happen and so there exists a prime $P$ such that $\langle z' \rangle \varsubsetneq P \varsubsetneq \langle z-\alpha, z' \rangle $. Note that these three ideals contain $z'$, but not $z$. Hence, applying the above bijection $\varphi$ to this chain, and then contracting with the center of $A$, we obtain a chain $ \langle 0 \rangle_{Z(A)} \varsubsetneq Q=\varphi(P) \cap Z(A) \varsubsetneq \langle \pi(z)-\alpha \rangle_{Z(A)} $  of prime ideals in $Z(A)=\mathbb{C}[\pi(z)^{\pm 1} ]$. Naturally this is a contradiction and so we have proved that $\haut(\langle z-\alpha ,z' \rangle) =2$, as desired. \fin
\\$ $


We now investigate the $\hc$-stratum associated to $\langle 0 \rangle$. By using similar arguments, we obtain the following description of $\spec_{\langle 0 \rangle}(U^+)$.

\begin{prop} 
\label{stratezero}
Let $\mathcal{P}$ be the set of those unitary irreductible polynomials $P(z,z') \in \mathbb{C}[z,z'] $ with $P(z,z') \neq z$ and $P(z,z') \neq z'$. 
\\Then $\spec_{\langle 0 \rangle } (U^+) = \{ \langle 0 \rangle  \} \cup \{ \langle P(z,z') \rangle \mid P(z,z') \in \mathcal{P} \} \cup \{ \langle z-\alpha ,z'-\beta \rangle \mid \alpha,\beta \in \mathbb{C}^* \} $.
\\Further, for all $P(z,z') \in \mathcal{P}$, $\haut(\langle P(z,z') \rangle )=1$, and, for all $\alpha,\beta \in \mathbb{C}^*$, $\haut(\langle z-\alpha ,z'-\beta \rangle)=2$.
\end{prop}
\preuve The proof of this result is similar to the proof of Proposition \ref{stratezprime}. Nevertheless, we include the proof since this proposition will play a crucial role in the sequel of this paper.

Observe first that the prime ideals in $\spec_{\langle 0 \rangle } (U^+)$ do not contain $z$ and $z'$. Indeed, assume that this is not the case, that is, assume that there exists $Q \in \spec_{\langle 0 \rangle } (U^+)$ with $zz' \in Q$. Then, since $zz'$ is an $\hc$-eigenvector, we have $zz' \in \bigcap_{h\in \hc} h.Q = \langle 0 \rangle$. This is a contradiction and so we have just proved that $\spec_{\langle 0 \rangle } (U^+) \subseteq \{ Q \in \spec(U^+) \mid z,z' \notin Q \}$. On the other hand, if $Q$ is a prime ideal of $U^+$ such that $z$, $z' \notin Q$, then $\bigcap_{h\in \hc} h.Q $ is an $\hc$-invariant prime ideal of $U^+$ that does not contain $z$ nor $z'$. In view of the list of $\hc$-invariant prime ideals of $U^+$ (see Proposition \ref{hpremier}), the only possibility is $\bigcap_{h\in \hc} h.Q =\langle 0 \rangle $, so that $Q \in \spec_{\langle 0 \rangle } (U^+)$. To resume, we have shown that  $\spec_{\langle 0 \rangle } (U^+) = \{ Q \in \spec(U^+) \mid z,z' \notin Q \}$.

Since $z,z'$ belong to the center of $U^+$, $\{z^iz^{'j} \mid i,j \in \mathbb{N} \}$ is a right denominator set in $U^+$; the corresponding localisation of $U^+$ will be denoted $\displaystyle{A:=U^+[z^{-1},z^{'-1}] }$. Note that, since $z$ and $z'$ are $\hc$-eigenvectors, the torus $\hc$ still acts rationally by automorphisms on $A$ (see \cite[Exercise II.3.A]{bg}). Moreover, it follows from the previous study (and from classical results of non-commutative localisation theory) that the map $\varphi : Q \rightarrow Q[z^{-1},z^{'-1}] $ is an increasing bijection from $\spec_{\langle 0 \rangle } (U^+)$ onto $\spec \left(A \right)$.

Before describing the prime spectrum of $A$, we establish that the algebra $A$ is $\hc$-simple in the sense of \cite[II.1.8]{bg}, that is, we show that $A$ has only one (two-sided) proper $\hc$-invariant ideal: $\langle 0 \rangle_A$. First, since $z$ and $z'$ are $\hc$-eigenvectors, we deduce from \cite[Exercise II.1.J]{bg} that the torus $\hc$ still acts by automorphisms on $A$ and that the bijection $\varphi$ induces a bijection between the set of those $\hc$-invariant prime ideals of $U^+$ that do not contain $z$ nor $z'$ and the set of $\hc$-invariant prime ideals of $A$. Since the set of those $\hc$-invariant prime ideals of $U^+$ that do not contain $z$ nor $z'$ is reduced to $\{ \langle 0 \rangle \}$ (see Proposition \ref{hpremier}), we obtain that $A$ has only one $\hc$-invariant prime ideal: $\langle 0 \rangle_A$. Now, every $\hc$-invariant proper ideal of $A$ is contained in an $\hc$-invariant prime ideal of $A$, so that $\langle 0 \rangle_A$ is the unique $\hc$-invariant proper ideal of $A$. In other words, $A$ is $\hc$-simple, as desired.

We are now able to describe the prime spectrum of $A$ and the $\hc$-stratum of $U^+$ associated to $\langle 0 \rangle$. Since the action of $\hc$ on $A$ is rational (see \cite[Exercise II.3.A]{bg}) and since $A$ is $\hc$-simple, it follows from  \cite[Corollary II.3.9]{bg} that extension and contraction provide mutually inverse bijections between $\spec(A)$ and $\spec(Z(A))$. Now, it follows from \cite[Lemma 3.1]{andrusdumas} that the center of $U^+$  is the polynomial algebra $Z(U^+)=\mathbb{C}[z,z']$, so that $Z(A)=\mathbb{C}[z^{\pm 1},z^{'\pm 1}]$. Since $\mathbb{C}$ is algebraically closed, we obtain that $\spec (Z(A))= \{ \langle 0 \rangle  \} \cup \{ \langle P(z,z') \rangle_{Z(A)} \mid P(z,z') \in \mathcal{P} \} \cup \{ \langle z-\alpha ,z'-\beta \rangle_{Z(A)} \mid \alpha,\beta \in \mathbb{C}^* \} $ and so we deduce from the previous study that
$\spec(A) =  \{ \langle 0 \rangle_A  \} \cup \{ \langle P(z,z') \rangle_A \mid P(z,z') \in \mathcal{P} \} \cup \{ \langle z-\alpha ,z'-\beta \rangle_A  \mid \alpha,\beta \in \mathbb{C}^* \} $, and next that $\spec_{\langle 0 \rangle } (U^+) =  \{ \langle 0 \rangle  \} \cup \{ \langle P(z,z') \rangle_A \cap U^+ \mid P(z,z') \in \mathcal{P} \} \cup \{ \langle z-\alpha ,z'-\beta \rangle_A \cap U^+ \mid \alpha,\beta \in \mathbb{C}^* \} $. Naturally, $\langle 0 \rangle_A \cap U^+=\langle 0 \rangle$. Thus, to achieve the proof of the first part of the proposition, it just remains to show that $\langle P(z,z') \rangle_A \cap U^+ = \langle P(z,z') \rangle $ for all $P(z,z') \in \mathcal{P}$, and  $\langle z-\alpha ,z'-\beta \rangle_A \cap U^+ = \langle z-\alpha ,z'-\beta \rangle$ for all $ \alpha,\beta \in \mathbb{C}^*$. This is what we do now.

Fix $P(z,z') \in \mathcal{P}$. We show that $ \langle P(z,z')  \rangle_A  \cap U^+ = \langle P(z,z') \rangle$. Clearly we have $ \langle P(z,z')  \rangle_A  \cap U^+ \supseteq \langle P(z,z') \rangle$  and so we just have to establish the reverse inclusion. Let $x \in \langle P(z,z')  \rangle_A \cap U^+$. Then there exist $s,t \in \mathbb{N}$ and $a \in U^+$ such that $x z^s z^{'t} =a P(z,z')$. Choose $(s,t)$ minimal (for the lexicographic order on $\mathbb{N}^2$). If, for instance, $s > 0$, then 
$a P(z,z') \in \langle z \rangle$. Since this ideal is completely prime, we get $a \in \langle z \rangle$ or $P(z,z') \in \langle z \rangle$. Now, since $P(z,z') \in \mathcal{P}$, the last case cannot happen, so that $a \in \langle z \rangle$. Thus, since $z$ is central, there exists $b \in U^+$ such that $a=bz$ and so $x z^s z^{'t}=aP(z,z')=bzP(z,z')$. Hence, we have $x z^{s-1} z^{'t}=bP(z,z')$. This contradicts the minimality of $(s,t)$. Thus $s=t=0$ and so $x =aP(z,z')$ belongs to the ideal of $U^+$ generated by $P(z,z')$. So we just prove that $ \langle P(z,z')  \rangle_A  \cap U^+ = \langle P(z,z') \rangle$.

Fix now $\alpha,\beta \in \mathbb{C}^*$.  Clearly we have $ \langle z-\alpha ,z'-\beta \rangle_A  \cap U^+ \supseteq \langle z-\alpha ,z'-\beta \rangle$ and so, once again, we just need to establish the reverse inclusion. Let $x \in \langle z-\alpha ,z'-\beta \rangle_A \cap U^+$. Then there exist $s,t \in \mathbb{N}$ and $a,b \in U^+$ such that $x z^s z^{'t} =a (z-\alpha) +b (z'-\beta)$. Choose $(s,t)$ minimal (for the lexicographic order on $\mathbb{N}^2$). If, for instance, $s > 0$, then, in $U^+/ \langle z' - \beta \rangle$ which is a domain (because of the previous point), we have:
$$\psi(x) \psi(z)^s \psi(z')^{t}  = \psi(a) (\psi(z)-\alpha),$$
where $\psi$ denotes the canonical surjection from $U^+$ onto $U^+ / \langle z' - \beta \rangle$. Since $s>0$ and $\alpha \neq 0$, we get $\psi(a) \in \langle \psi(z) \rangle$. Now, since $\psi(z)$ is a central element of $U^+/ \langle z' - \beta \rangle$, there exists $a' \in U^+$ such that $\psi(a)=\psi (a') \psi(z)$. 
Then, $\psi(x) \psi(z)^s \psi(z')^{t} = \psi(a) (\psi(z)-\alpha)= \psi(a') \psi(z) (\psi(z)-\alpha)$. Since $\psi(z)$ is a non-zero element of the domain $U^+/ \langle z' - \beta \rangle$, we obtain that $\psi(x) \psi(z)^{s-1} \psi(z')^{t} =  \psi(a') (\psi(z)-\alpha)$. Hence there exists $b' \in U^+$ such that  
$x z^{s-1} z^{'t} =a' (z-\alpha) +b' (z'-\beta)$. This contradicts the minimality of $(s,t)$. Thus, we have $s=t=0$, and so $x =a (z-\alpha) +b (z'-\beta) \in \langle z-\alpha ,z'-\beta \rangle$, as desired. This achieves the first part of the proposition.

The second part of the proposition can be easily obtained by using the ordered bijection $\phi : \spec_{\langle 0 \rangle } (U^+) \rightarrow \spec(Z(A))=\spec(\mathbb{C}[z^{\pm 1}, z^{'\pm 1}])$ defined by 
$\phi(Q) = Q [z^{-1},z^{'-1}] \cap Z(A)$ for all $Q \in \spec_{\langle 0 \rangle } (U^+)$. \fin


\subsubsection{Height $3$ prime ideals of $U^+$.}
$ $

 Using the above description of the 8 $\hc$-strata of $\spec(U^+)$, we easily obtain the following result.

\begin{prop}
\label{haut3}
$U^+$ has only two height $3$ prime ideals: $\langle e_1 \rangle$ et $\langle e_2 \rangle$. 
\\In particular, $z$ and $z'$ belong to every prime ideal of height greater than or equal to $3$.
\end{prop}
\preuve It follows from Proposition \ref{stratezprime} and \ref{stratezero} that que the prime ideals of $U^+$ that do not contain $z$ have height at most $2$. Hence $z$ belongs to every height $3$ prime ideal of $U^+$. Thus, the height $3$ prime ideals of $U^+$ are exactly the inverse images by the canonical surjection $\pi : U^+ \rightarrow U^+/\langle z \rangle$ of the height $2$ prime ideals of $U^+/\langle z \rangle$. Since $U^+/\langle z \rangle$ is isomorphic to the quantum Heisenberg algebra $\mathbb{H}$ (see Section \ref{paraheisenberg}), we deduce from \cite[Th\'eor\`eme 2.4]{malliavin} that $U^+/\langle z \rangle$ has only two height $2$ prime ideals: 
$\langle \pi(e_1) \rangle$ and $\langle \pi(e_2) \rangle$. Hence, $U^+$ has only two height $3$ prime ideals : $\langle e_1 \rangle$ et $\langle e_2 \rangle$. \fin
\\$ $


\subsection{(Left) primitive and maximal ideals of $U^+$.}
$ $

We have seen that the $\hc$-strata $\spec_J(U^+)$ ($J \in \hc$-$\spec(U^+)$) form a partition of $\spec(U^+)$:
$$\spec(U^+) = \bigsqcup_{J \in \hc \mbox{-}\spec(U^+)} \spec_J(U^+).$$
Naturally, this partition induces a partition of the set $\prim(U^+)$ of all (left) primitive ideals of $U^+$ as follows. 
For all $J \in \hc$-$\spec(U^+)$, we set $\prim_J(U^+):=\spec_J(U^+) \cap \prim(U^+)$. Then it is obvious that 
the $\hc$-strata $\prim_J(U^+)$ ($J \in \hc$-$\spec(U^+)$) form a partition of $\prim(U^+)$:
$$\prim(U^+) = \bigsqcup_{J \in \hc \mbox{-}\spec(U^+)} \prim_J(U^+).$$
We will now make precise this partition.
\\$ $

Since $\mathbb{C}$ is uncountable and since the Noetherian domain $U^+$ is generated (as an algebra) by a finite number of elements, it follows from \cite[Proposition II.7.16]{bg} that the algebra $U^+$ satisfies the Nullstellensatz over $\mathbb{C}$ (see \cite[II.7.14]{bg}). Further the set of $\hc$-invariant prime ideals of $U^+$ is finite. Thus we deduce from \cite[Theorem II.8.4]{bg} that $\prim_J(U^+)$ $(J \in \hc$-$\spec(U^+)$) coincides with those primes 
in $\spec_J(U^+)$ that are maximal in $\spec_J(U^+)$. Hence, we deduce from the above description of the $\hc$-strata of $\spec(U^+)$ (see Proposition \ref{stratez}, \ref{stratezprime} and \ref{stratezero}) the following description of the $\hc$-strata of $\prim(U^+)$.
\\$ $

\begin{prop}
\label{prim}
$$\begin{array}{lll}
\prim_{\langle 0 \rangle } (U^+) & = & \{ \langle z-\alpha ,z'-\beta \rangle \mid \alpha,\beta \in \mathbb{C}^* \} ,\\
\prim_{\langle z' \rangle } (U^+) & = &  \{ \langle z-\alpha ,z' \rangle \mid \alpha \in \mathbb{C}^* \} ,\\
\prim_{\langle z \rangle } (U^+) & = & \{ \langle z ,z'-\beta \rangle \mid \beta \in \mathbb{C}^* \} ,\\
\prim_{\langle e_3 \rangle } (U^+)  & = & \{ \langle e_3 \rangle \} ,\\
\prim_{\langle \ov{e_3} \rangle } (U^+) & = & \{ \langle \ov{e_3} \rangle \},\\
\prim_{\langle e_1 \rangle } (U^+) &  = &  \{ \langle e_1,e_2-\beta \rangle \mid \beta \in \mathbb{C}^* \} ,\\
\prim_{\langle e_2 \rangle } (U^+) & = & \{ \langle e_1-\alpha , e_2 \rangle \mid \alpha \in \mathbb{C}^* \} ,\\
\prim_{\langle e_1,e_2 \rangle } (U^+) & = & \{ \langle e_1,e_2 \rangle \}.
\end{array}$$
\end{prop}

Let $\mx(U^+)$ denote the set of maximal ideals of $U^+$. Among the primitive ideals, only two are not maximal. Indeed,  we have:

\begin{prop}
\label{max}
$\mx(U^+) $ is the disjoint union of the following $\hc$-strata of $\prim(U^+)$: $\{ \langle z-\alpha ,z'-\beta \rangle \mid \alpha,\beta \in \mathbb{C}^* \}$, $\{ \langle z ,z'-\beta \rangle \mid \beta \in \mathbb{C}^* \} $, $\{ \langle z-\alpha ,z' \rangle \mid \alpha \in \mathbb{C}^* \} $, $\{ \langle e_1, e_2-\beta \rangle  \mid \beta \in \mathbb{C}^* \} $, $\{ \langle e_1 - \alpha , e_2 \rangle  \mid \alpha \in \mathbb{C}^* \} $ and $\{ \langle e_1, e_2 \rangle \}$.
\\In other words, in $U^+$, there are only two primitive ideals that are not maximal: $\langle e_3 \rangle$ and $\langle \ov{e_3} \rangle$.
\end{prop}  
\preuve First, $\langle e_3 \rangle$ and $\langle \ov{e_3} \rangle$ are not maximal since they are stricly contained  in the augmentation ideal $\langle e_1,e_2 \rangle$.

Next, $\langle e_1, e_2-\beta \rangle$ ($\beta \in \mathbb{C}^*$), $\langle e_1 - \alpha , e_2 \rangle$ ($ \alpha \in \mathbb{C}^*$) and $\langle e_1, e_2 \rangle$ are maximal since they have maximal height.

Now let $\alpha \in \mathbb{C}^*$. Recall that $\haut(\langle z-\alpha ,z' \rangle)=2$. If $\langle z-\alpha ,z' \rangle$ is not maximal, then there exists a height $3$ prime $P$ such that $\langle z-\alpha ,z' \rangle \subseteq P$. In particular, we have $z-\alpha \in P$. On the other hand, since every height $3$ prime contains $z$ (see Proposition \ref{haut3}), we also have $z \in P$. This is a contradiction since $\alpha \in \mathbb{C}^*$. 
Hence $\langle z-\alpha ,z' \rangle$ is maximal

We prove with the same arguments that the primitive ideals $\langle z-\alpha ,z'-\beta \rangle$ ($\alpha,\beta \in \mathbb{C}^*$) and $\langle z ,z'- \beta \rangle$ ($\beta \in \mathbb{C}^*$) are also maximal. \fin
\\$ $

We finish this section by computing the Gelfand-Kirillov dimension of the simple factor algebras of $U^+$. 
\\$ $

\begin{prop}
\label{gkdimmax}
$ $
\\$\bullet$ $ \displaystyle{ \gk \left( \frac{U^+}{\langle z-\alpha ,z'-\beta \rangle} \right) = 2 } $ for all $\alpha, \beta \in \mathbb{C}^*$.
\\$\bullet$ $ \displaystyle{ \gk \left( \frac{U^+}{\langle z ,z'-\beta \rangle} \right) = 2 }$ for all $\beta \in \mathbb{C}^*$.
\\$\bullet$ $ \displaystyle{ \gk \left( \frac{U^+}{\langle z-\alpha ,z' \rangle} \right) = 2 }$ for all $\alpha \in \mathbb{C}^*$.
\\$\bullet$ $ \displaystyle{ \gk \left( \frac{U^+}{\langle e_1 - \alpha , e_2 \rangle } \right) = 0 }$ for all $\alpha, \beta \in \mathbb{C}^*$.
\\$\bullet$ $ \displaystyle{ \gk \left( \frac{U^+}{\langle e_1  , e_2 -\beta \rangle } \right) = 0 }$ for all $\beta \in \mathbb{C}^*$.
\\$\bullet$ $ \displaystyle{ \gk \left( \frac{U^+}{\langle e_1  , e_2 \rangle } \right) = 0 }$.
\end{prop}
\preuve From \cite[Theorem 4.8]{goodlenCAT}, Tauvel's height formula holds in $U^+$. So, for all prime $P$ in $U^+$, we have 
$$\gk \left( \frac{U^+}{P} \right) + \haut(P) = \gk (U^+)=4. $$
In Proposition \ref{stratez}, \ref{stratezprime} and \ref{stratezero}, we have calculated the height of all primes of $U^+$ and so the result easily follows from 
these three Proposition \ref{stratez}, \ref{stratezprime} and \ref{stratezero}. \fin

\begin{rem}
Except for $\langle z,z' \rangle_{Z(U^+)}$, the maximal ideals of the center $Z(U^+)=\mathbb{C}[z,z']$ of $U^+$ extend to height $2$ maximal ideals of $U^+$. However they do not appear in the same $\hc$-strata. In particular, the height $2$ maximal ideals that contain $z$ and those that contain $z'$ are not in the same $\hc$-strata. This indicates that $z$ and $z'$ must be distinguished: they do not have the same status in $U^+$. This observation, that will play a crucial role in the calculation of the automorphism group of $U^+$, can be explained as follows: the two indeterminates $z$ and $z'$ that generate the center do not play the same role since $z$ comes directly from the Lie algebra of type $B_2$ contrary to $z'$ which only appears at the (quantized) enveloping algebra level.
\end{rem}  

\section{Simple factor algebras of Gelfand-Kirillov dimension $2$ of $U^+$.}
$ $

This paragraph is devoted to the simple factors of Gelfand-Kirillov dimension $2$ of $U^+$, that is, in view of Proposition \ref{max} and \ref{gkdimmax}, the factor algebras $U^+ / \langle z-\alpha ,z'-\beta \rangle$ ($(\alpha , \beta) \in \mathbb{C}^2 \setminus \{(0,0) \}$). Recall that, in the classical case, the simple factors of Gelfand-Kirillov dimension $2$ of the enveloping algebra of a Lie algebra of type $B_2$ are isomorphic to the first Weyl algebra $A_1(\mathbb{C})$ (see \cite{Dixmier}). Hence, our aim in this section is to compare the properties 
of the algebras $U^+ / \langle z-\alpha ,z'-\beta \rangle$ ($(\alpha ,
\beta) \in \mathbb{C}^2 \setminus \{(0,0) \}$) with those of the first
Weyl algebra $A_1(\mathbb{C})$. In particular, we will prove that the
centers of such simple factor algebras are reduced to $\mathbb{C}$ and
we will calculate the groups of units of these algebras. This study
will suggest that we distinguish between three families of such
algebras:\\$ $
\\$\bullet$ The algebras $U^+ / \langle z ,z'-\beta \rangle$ with
$\beta \neq 0$  are isomorphic to simple factor algebras of
Gelfand-Kirillov dimension $2$ of the quantum Heisenberg
algebra. These algebras are the so-called Weyl-Hayashi algebras (see
\cite{hayashi}) and have been studied from a ring-theoritical point of
view by Alev and Dumas (see \cite{alevdumasrigidite}), Kirkman and
Small (see \cite{kirkman}), and  Malliavin (see \cite{malliavin}).
They are simple, but they contain invertible elements that are not
scalars. In fact, they can be presented as generalized Weyl algebras
(GWA for short) over a Laurent polynomial ring in one variable (see
Section \ref{GWA} for the definition of a GWA).
\\$ $
\\$\bullet$ The algebras $U^+ / \langle z-\alpha ,z' \rangle$ with
$\alpha \neq 0$ are also simple and contain invertible elements that
are not scalars. In fact, they can also be presented as generalized
Weyl algebras over a Laurent polynomial ring in one
variable. However, they are not isomorphic to any algebra in
the first family. 
\\$ $
\\$\bullet$ The algebras $U^+ / \langle z-\alpha ,z'-\beta \rangle$
with $\alpha, \beta \neq 0$ are simple, their centers are reduced to
$\mathbb{C}$ and they do not contain non-trivial invertible
elements. Further, they cannot be presented as generalized Weyl
algebras over a polynomial ring in one variable nor a Laurent
polynomial ring in one variable. Algebras in this third family provide
good quantizations of the first Weyl algebra.


\subsection{Basis of $U^+ / \langle z-\alpha ,z'-\beta \rangle$ ($(\alpha , \beta) \in \mathbb{C}^2 \setminus \{(0,0)\} $).}
\label{soussectionbase} 
$ $

Fix $(\alpha , \beta) \in \mathbb{C}^2 \setminus \{(0,0)\}$. We set 
$B_{\alpha}:= U^+ / \langle z- \alpha \rangle$ and $A_{\alpha,\beta} :=U^+ / \langle z-\alpha ,z'-\beta \rangle  $; these algebras are Noetherian domains. If $x \in U^+$, we denote by $\widehat{x}$ (resp. $\widetilde{x}$) the canonical image of $x$ in $B_{\alpha}$ (resp. $A_{\alpha,\beta}$). Naturally, in the sequel, we will identify $A_{\alpha,\beta}$ with $B_{\alpha}/ \langle \widehat{z'}- \beta \rangle $.

Recall (see \ref{quotientzalpha}) that the family $(\widehat{e_3}^k\widehat{e_1}^i\widehat{e_2}^j)_{i,j,k \in \mathbb{N}}$ is a basis of $B_{\alpha}$ as $\mathbb{C}$-vector space.

On the other hand, we deduce from the commutation relations between the generators $e_1,e_2,e_3$ of $U^+$ and $z \in Z(U^+)$ that $A_{\alpha,\beta}$ is generated as an algebra by $\teun,\tedeux $ and $\tetrois$, and that, in $A_{\alpha,\beta}$, we have the following commutation relations.
$$ \teun \tetrois = q^{-2} \tetrois \teun \mbox{, } \hspace{1cm} \tedeux \tetrois = q^2 \tetrois \tedeux + \alpha 
\mbox{, }\hspace{1cm} \tedeux \teun = q^{-2} \teun \tedeux - q^{-2}\tetrois $$
$$\mbox{ and }\hspace{1cm} \tetrois^2+c_1 \tetrois \teun \tedeux +\alpha c_2 \teun + \beta c_3 =0,$$
where $c_1= q^4-1$, $c_2=q^2(q^2+1)$ and $c_3=\frac{q^6}{1-q^2}$.
\\$ $

As observed by Kirkman and Small (see \cite{kirkman}), the Weyl-Hayashi algebras $A_{0,\beta}$ ($\beta \in \mathbb{C}$) does not become the first Weyl algebra $A_1(\mathbb{C})$ when $q=1$. One of the advantage of considering the algebras 
$A_{\alpha,\beta}$ ($(\alpha,\beta) \in \mathbb{C}^2 \setminus \{(0,0)\}$) is that for suitable values of the parameters $\alpha$ and $\beta$ 
the algebra $A_{\alpha,\beta}$ becomes the first Weyl algebra when $q=1$. More precisely, one can check that

\begin{obs}
$A_{1,0}$, $A_{1,(q-1)^2}$... become the first Weyl algebra $A_1(\mathbb{C})$ when $q=1$.
\end{obs}
$ $

The following result provides a linear basis of $A_{\alpha,\beta}$.

\begin{prop}
\label{base}
Set $\mathcal{F}:= (\teun^i \tedeux^j)_{i,j \in \mathbb{N}} \cup (\tetrois \teun^i \tedeux^j)_{i,j \in \mathbb{N}}$. 
\\$\mathcal{F}$ is a basis of the $\mathbb{C}$-vector space $A_{\alpha,\beta}$.
\end{prop}
\preuve Fix $(\alpha , \beta) \in \mathbb{C}^2\setminus \{(0,0)\}$. We begin by showing that the $\mathbb{C}$-vector space  $A_{\alpha,\beta}$ is generated by $\mathcal{F}$. First, since $A_{\alpha,\beta} \simeq B_{\alpha}/\langle \widehat{z'}-\beta \rangle$ and since the family $(\widehat{e_3}^k\widehat{e_1}^i\widehat{e_2}^j)_{i,j,k \in \mathbb{N}}$ is a linear basis of $B_{\alpha}$ (see Section \ref{quotientzalpha}), the family $(\tetrois^k \teun^i\tedeux^j)_{i,j,k \in \mathbb{N}}$ is a generating set of $A_{\alpha,\beta}$. So it is sufficient to establish that, for all $k,i,j \in \mathbb{N}$, 
$\tetrois^k \teun^i\tedeux^j$ can be written as a (finite) linear combination over $\mathbb{C}$ of elements of $\mathcal{F}$. To do this, we proceed by induction on $k$. 

If $k=0$ or $1$, it is obvious. We now assume that $k \geq 2$. It follows from the commutation relations in $A_{\alpha,\beta}$ that we have:
\begin{eqnarray*}
\tetrois^k \teun^i\tedeux^j & = & q ^{4i} \tetrois^{k-2} \teun^i \tetrois^{2} \tedeux^j \\
 & = & q ^{4i} c_1 \tetrois^{k-2} \teun^i \tetrois \teun \tedeux^{j+1} + q ^{4i}c_2 \alpha \tetrois^{k-2} \teun^{i+1}  \tedeux^j + q ^{4i} c_3 \beta \tetrois^{k-2} \teun^i  \tedeux^j \\
  & = & q ^{2i} c_1 \tetrois^{k-1} \teun^{i+1}  \tedeux^{j+1} + q ^{4i}c_2 \alpha \tetrois^{k-2} \teun^{i+1}  \tedeux^j + q ^{4i} c_3 \beta \tetrois^{k-2} \teun^i  \tedeux^j \\
\end{eqnarray*}
Then we deduce from the induction hypothesis that $\tetrois^k \teun^i\tedeux^j$ can be written as a (finite) linear combination over $\mathbb{C}$ of elements of $\mathcal{F}$ as required. This achieves the induction and so we have just proved that $\mathcal{F}$ is a generating family of $A_{\alpha,\beta}$ viewed as a $\mathbb{C}$-vector space.
\\$ $

We now establish that the family $\mathcal{F}$ is linearly independent. 
Let $I$ and $I'$ be two finite subsets of $\mathbb{N}^2$, and $(\alpha_{i,j})_{(i,j) \in I}$, $(\beta_{i,j})_{(i,j) \in I'}$ two families of complex numbers. We assume that
$$\sum_{(i,j) \in I} \alpha_{i,j} \teun^i \tedeux^j + \sum_{(i,j) \in I'} \beta_{i,j} \tetrois \teun^i \tedeux^j=0.$$
Then, since the family $(\widehat{e_3}^k\widehat{e_1}^i\widehat{e_2}^j)_{i,j,k \in \mathbb{N}}$ is a linear basis of $B_{\alpha}$, 
and since $A_{\alpha,\beta} \simeq B_{\alpha}/ \langle \widehat{z'}- \beta \rangle =B_{\alpha}/ \langle \hetrois^2+c_1 \hetrois \heun \hedeux +\alpha c_2 \heun + \beta c_3 \rangle$, there exist a finite subset $\Gamma$ of $\mathbb{N}^3$ and a family 
  $(\gamma_{k,i,j})_{(k,i,j) \in \Gamma}$ of complex numbers such that, in $B_{\alpha}$, we have:
\begin{eqnarray}
\label{egaliteliberte}
\sum_{(i,j) \in I} \alpha_{i,j} \heun^i \hedeux^j + \sum_{(i,j) \in I'} \beta_{i,j} \hetrois \heun^i \hedeux^j 
=\left(\hetrois^2+c_1 \hetrois \heun \hedeux +\alpha c_2 \heun + \beta c_3 \right) \sum_{(k,i,j) \in \Gamma} \gamma_{k,i,j} \hetrois^k \heun^i \hedeux^j .
\end{eqnarray}
Further the commutation relations of Lemma \ref{relationsdansU} show that, in $B_{\alpha}$, we have  
\begin{enumerate}
\item $\displaystyle{\hetrois \heun^k  =  q^{2k} \heun^k \hetrois}$ and $\displaystyle{ \heun \hetrois^k  =  q^{-2k} \hetrois^k \heun}$.
\item $\displaystyle{\hedeux \hetrois^k  =  q^{2k} \hetrois^k \hedeux + \alpha \frac{q^{2k}-1}{q^2-1}\hetrois^{k-1}}$.
\item $\displaystyle{\hedeux \heun^k  =  q^{-2k} \heun^k \hedeux - q^{-2} \frac{q^{-4k}-1}{q^{-4}-1} \hetrois \heun^{k-1}}$,
\end{enumerate}
for all $k \in \mathbb{N}^*$. We easily deduce from the above relations that
$$\hetrois \heun \hedeux \hetrois^k \heun^i \hedeux^j 
= q^{-2i} \hetrois^{k+1} \heun^{i+1} \hedeux^{j+1} 
-q^{-4} \frac{q^{-4i}-1}{q^{-4}-1} \hetrois^{k+2} \heun^{i} \hedeux^{j} 
+ \alpha q^{-2(k-1)} \frac{q^{2k}-1}{q^{2}-1} \hetrois^{k} \heun^{i+1} \hedeux^{j},$$
and next that
\begin{eqnarray*}
\sum_{(i,j) \in I} \alpha_{i,j} \heun^i \hedeux^j  +   \sum_{(i,j) \in I'} \beta_{i,j} \hetrois \heun^i \hedeux^j 
 & = &  \sum_{(k,i,j) \in \Gamma} \gamma_{k,i,j} \hetrois^{k+2} \heun^i \hedeux^j \\
& + &  \sum_{(k,i,j) \in \Gamma} q^{-2i} c_1 \gamma_{k,i,j} \hetrois^{k+1} \heun^{i+1} \hedeux^{j+1} \\
& - &  \sum_{(k,i,j) \in \Gamma} q^{-4} \frac{q^{-4i}-1}{q^{-4}-1} c_1 \gamma_{k,i,j} \hetrois^{k+2} \heun^i \hedeux^j \\
& + &  \sum_{(k,i,j) \in \Gamma} \alpha q^{-2(k-1)} \frac{q^{2k}-1}{q^{2}-1}  c_1 \gamma_{k,i,j} \hetrois^{k} \heun^{i+1} \hedeux^j \\
& + &  \sum_{(k,i,j) \in \Gamma} \alpha c_2 q^{-2k} \gamma_{k,i,j} \hetrois^{k} \heun^{i+1} \hedeux^j \\
& + &  \sum_{(k,i,j) \in \Gamma} \beta c_3 \gamma_{k,i,j} \hetrois^{k} \heun^i \hedeux^j \\
\end{eqnarray*}
Assume that there exists $(k,i,j) \in \Gamma$ such that $\gamma_{k,i,j} \neq 0$. Then we denote by $(w,u,v)$ the greatest element (relative to the lexicographic order on $\mathbb{N}^3$) of $\Gamma$ such that $\gamma_{w,u,v} \neq 0$. Since the family $(\widehat{e_3}^k\widehat{e_1}^i\widehat{e_2}^j)_{i,j,k \in \mathbb{N}}$ is a linear basis of $B_{\alpha}$ (see Section \ref{quotientzalpha}), we can identify the coefficients of $\hetrois^{w+2} \heun^u \hedeux^v$ in the previous equality. This leads to
$$0=\gamma_{w,u,v} -q^{-4} \frac{q^{-4u}-1}{q^{-4}-1} c_1 \gamma_{w,u,v}. $$
Thus, since $\gamma_{w,u,v} \neq 0$, we have $1 -q^{-4} \frac{q^{-4u}-1}{q^{-4}-1} c_1 =0$. Further, $c_1 = q^4 -1$ and so we get $q^{-4u}=0$. This is a contradiction since $q$ is nonzero. Hence $\gamma_{k,i,j} = 0$ for all $(k,i,j) \in \Gamma$, and so it follows from the above equality (\ref{egaliteliberte}) that: 
$$\sum_{(i,j) \in I} \alpha_{i,j} \heun^i \hedeux^j + \sum_{(i,j) \in I'} \beta_{i,j} \hetrois \heun^i \hedeux^j =0.$$
One more time, since the family $(\widehat{e_3}^k\widehat{e_1}^i\widehat{e_2}^j)_{i,j,k \in \mathbb{N}}$ is a linear basis of $B_{\alpha}$, this implies that $\alpha_{i,j} = 0$ for all $(i,j) \in I$ and $\beta_{i,j} = 0$ for all $(i,j) \in I'$, so that $\mathcal{F}$ is linearly independent, as desired. \fin 


\subsection{The group of invertible elements of $U^+ / \langle z-\alpha ,z'-\beta \rangle$ ($(\alpha , \beta) \in \mathbb{C}^2 \setminus \{(0,0) \}$).}
$ $

 Fix $(\alpha , \beta) \in \mathbb{C}^2 \setminus \{(0,0) \}$. We keep the conventions and notations introduced in the previous section. In particular, we still set $B_{\alpha}:= U^+ / \langle z- \alpha \rangle$, $A_{\alpha,\beta} :=U^+ / \langle z-\alpha ,z'-\beta \rangle  $ and $\widehat{x}$ (resp. $\widetilde{x}$) the canonical image of $x$ in $B_{\alpha}$ (resp. $A_{\alpha,\beta}$). 

Before describing the group of invertible elements of $A_{\alpha,\beta}$, we show that a suitable localisation of this algebra is isomorphic to a quantum torus. More precisely, we denote by $ \mathbb{C}_{q^2}[x, y]$ the quantum plane, that is, the algebra generated by two indeterminates $x,y$ with the relation $xy=q^2 yx$, and we denote by $ \mathbb{C}_{q^2}[x^{\pm 1}, y^{\pm 1}]$ the quantum torus associated to the quantum plane $ \mathbb{C}_{q^2}[x, y]$, that is, the localisation of the quantum plane $ \mathbb{C}_{q^2}[x, y]$ by the multiplicative system generated by the normal elements $x$ and $y$. Then we have the following statement.


\begin{lem}
\label{propinter}
Set $\Sigma:=\{\lambda \tetrois^i \teun^j \mid \lambda \in \mathbb{C}^* \mbox{, } i,j \in \mathbb{N} \}$. 
\\$\Sigma$ is a  multiplicative system of regular elements of $A_{\alpha,\beta}$, that satisfies the Ore condition (both left and right) in $A_{\alpha,\beta}$.
\\Further, the map $g:\mathbb{C}_{q^2}[x^{\pm 1}, y^{\pm 1}] \rightarrow A_{\alpha,\beta}\Sigma^{-1}$ defined by 
$g(x)=\tetrois$ and $g(y) = \teun$ is an algebra isomorphism. 
\end{lem}  
\preuve Since $A_{\alpha,\beta}$ is a domain and since $\teun$ and $\tetrois$ are non-zero, $\Sigma$ is a  multiplicative system of regular elements of $A_{\alpha,\beta}$. Moreover, since the multiplicative system 
$\{\lambda e_3^i e_1^j \mid \lambda \in \mathbb{C}^* \mbox{, } i,j \in \mathbb{N} \}$ satisfies the Ore condition (both left and right) in $U^+$ (see \cite[3.1.4]{andrusdumas}), it is obvious that $\Sigma$ also satisfies the Ore condition (both left and right) in $A_{\alpha,\beta}$. This establishes the first part of the lemma.

Now, $\teun$ and $\tetrois$ are invertible in the Ore localisation $A_{\alpha,\beta}\Sigma^{-1}$. Further, we have: 
$\tetrois \teun = q^2 \teun \tetrois$. Hence, there exists an algebra homomorphism $g:\mathbb{C}_{q^2}[x^{\pm 1}, y^{\pm 1}] \rightarrow A_{\alpha,\beta}\Sigma^{-1}$ such that $g(x)=\tetrois$ and $g(y) = \teun$. It remains to prove that $g$ is actually an isomorphism. 

We first show that $g$ is onto. Recall (see \ref{soussectionbase}) that $A_{\alpha,\beta}$ is generated as an algebra  by $\teun$, $\tedeux$ and $\tetrois$, so that $A_{\alpha,\beta}\Sigma^{-1}$ is generated as an algebra by $\teun^{\pm 1}$, $\tedeux$ and $\tetrois^{\pm 1}$. Thus, in order to prove that $g$ is onto, it is sufficient to show that $\tedeux$ belongs to the sub-algebra of $A_{\alpha,\beta}\Sigma^{-1}$ generated by $\teun^{\pm 1}$ and $\tetrois^{\pm 1}$. This is what we do now.

Recall (see \ref{soussectionbase}) that, in $A_{\alpha,\beta}$ , 
$\teun \tetrois = q^{-2} \tetrois \teun$ and $ \tetrois^2+c_1 \tetrois \teun \tedeux +\alpha c_2 \teun + \beta c_3 =0$, where $c_1= q^4-1$, $c_2=q^2(q^2+1)$ and $c_3=\frac{q^6}{1-q^2}$.
We deduce from these relations that, in $A_{\alpha,\beta}\Sigma^{-1}$, we have:
$$\tedeux = -\frac{1}{c_1} \left(  q² \tetrois \teun^{-1} + q^2 \alpha c_2 \tetrois^{-1}  + q^{-2} \beta c_3 \tetrois^{-1} \teun^{-1} \right).$$
Thus, $\tedeux$ belongs to  the sub-algebra of $A_{\alpha,\beta}\Sigma^{-1}$ generated by $\teun^{\pm 1}$ and $\tetrois^{\pm 1}$, as required, and so $g$ is onto.

Now, to prove that the algebra homomorphism $g$ is actually an isomorphism, it is sufficient to show that the Noetherian domains $\mathbb{C}_{q^2}[x^{\pm 1}, y^{\pm 1}]$ and $A_{\alpha,\beta}\Sigma^{-1}$ have the same Gelfand-Kirillov dimension. Since it is well known that $\gk \left( \mathbb{C}_{q^2}[x^{\pm 1}, y^{\pm 1}] \right)=2$, it just remains to prove that the Gelfand-Kirillov dimension of $ A_{\alpha,\beta}\Sigma^{-1}$ is also 2. First, since $g$ is onto, we have :
\begin{eqnarray*}
\label{gk1}
\gk \left( A_{\alpha,\beta}\Sigma^{-1} \right) \leq \gk \left( \mathbb{C}_{q^2}[x^{\pm 1}, y^{\pm 1}] \right)=2.
\end{eqnarray*}
On the other hand, since $\gk \left( A_{\alpha,\beta} \right) = 2$, we have
\begin{eqnarray*}
\label{gk2}
\gk \left( A_{\alpha,\beta}\Sigma^{-1} \right) \geq 2.
\end{eqnarray*}
So $$  \gk \left( A_{\alpha,\beta}\Sigma^{-1} \right) =2 =\gk \left( \mathbb{C}_{q^2}[x^{\pm 1}, y^{\pm 1}] \right)$$
and $g$ is an algebra isomorphism, as desired. \fin \\$ $

An immediate consequence of Lemma \ref{propinter} is the following statement.

\begin{cor}
\label{corps}
The field of fractions of $A_{\alpha,\beta}$ ($(\alpha,\beta) \neq (0,0)$) is isomorphic to the field of fractions of the quantum plane  $\mathbb{C}_{q^2}[x,y]$.
\end{cor}
Before computing the group of invertible elements of
$A_{\alpha,\beta}$, we need some informations on the expression of
powers of $\tetrois$ in the linear basis of $A_{\alpha,\beta}$ that we
have described in Proposition \ref{base}.

\begin{lem}
\label{e3base}
Let $p$ be a positive integer. 
\begin{enumerate}
\item If $\beta \neq 0$, then we have 
$$\tetrois^{2p}=(-\beta c_3)^p + \sum_{i \geq 1, j \geq 0 }
  \alpha_{i,j} \teun^i \tedeux^j + \sum_{i \geq 1, j \geq 0}
  \beta_{i,j} \tetrois \teun^i \tedeux^j$$
and 
$$\tetrois^{2p+1}=(-\beta c_3)^p\tetrois + \sum_{i \geq 1, j \geq 0 }
  \alpha'_{i,j} \teun^i \tedeux^j + \sum_{i \geq 1, j \geq 0}
  \beta'_{i,j} \tetrois \teun^i \tedeux^j, $$
where $\alpha_{i,j}, \beta_{i,j},\alpha'_{i,j}, \beta'_{i,j} \in \mathbb{C}$ are equal to zero
  except for a finite number of them.
\item If $\beta=0$, then we have 
$$\tetrois^{2p}=  \sum_{i \geq p, j \geq 0 }
  \alpha_{i,j} \teun^i \tedeux^j + \sum_{i \geq p, j \geq 0}
  \beta_{i,j} \tetrois \teun^i \tedeux^j$$
and 
$$\tetrois^{2p+1}= \sum_{i \geq p, j \geq 0 }
  \alpha'_{i,j} \teun^i \tedeux^j + \sum_{i \geq p, j \geq 0}
  \beta'_{i,j} \tetrois \teun^i \tedeux^j, $$
where $\alpha_{i,j}, \beta_{i,j},\alpha'_{i,j}, \beta'_{i,j} \in \mathbb{C}$ are equal to zero
  except for a finite number of them, and $\alpha_{p,0} \neq 0$,
  $\beta'_{p,0} \neq 0$.
\end{enumerate}
\end{lem}
\preuve This lemma easily follows from a straightforward
induction. \fin \\$ $


In the sequel, if $R$ is an algebra, we denote by $U(R)$ the group of invertible elements of $R$. 
Since $U(\mathbb{C}_{q^2}[x^{\pm 1}, y^{\pm 1}]) = \{\lambda x^i y^j \mid \lambda \in \mathbb{C}^* \mbox{, } i,j \in \mathbb{Z} \}$, it follows from the above Lemma \ref{propinter} that $U(A_{\alpha,\beta}) \subseteq \{\lambda \tetrois^i \teun^j \mid \lambda \in \mathbb{C}^* \mbox{, } i,j \in \mathbb{Z} \}$. In fact, this inclusion is always strict.


\begin{prop}
\label{unite}
Let $(\alpha , \beta) \in \mathbb{C}^2 \setminus \{(0,0)\}$.
$$U(A_{\alpha,\beta}) = \left\{ \begin{array}{ll}
\{\lambda \tetrois^i \mid \lambda \in \mathbb{C}^* \mbox{, } i \in
\mathbb{Z} \} & \mbox{ if }\alpha=0 \\
\{\lambda (\teun^{-1} \tetrois^2)^i \mid \lambda \in \mathbb{C}^* \mbox{, } i \in \mathbb{Z} \} & \mbox{ if }\beta=0 \\
\mathbb{C}^* & \mbox{ otherwise.}
\end{array} \right.$$
\end{prop}
\preuve In the case where $\alpha=0$, $A_{\alpha,\beta}$ is a so-called
Weyl-Hayashi algebra whose group of invertible elements has been
computed for instance in \cite[Proof of Theorem
  2.5]{alevdumasrigidite}. Hence we assume from now that $\alpha \neq
0$. We proceed in three steps.\\$ $

\noindent $\bullet$ Step 1: we prove that $\teun \notin U(A_{\alpha,\beta})$ and  $\tetrois \notin 
U(A_{\alpha,\beta})$.\\$ $

Assume that $\teun \in U(A_{\alpha,\beta})$. Then there exists $a \in A_{\alpha,\beta}$ such that 
$\teun a =1$. It follows from Proposition \ref{base} that the family $\mathcal{F}:= (\teun^i \tedeux^j)_{i,j \in \mathbb{N}} \cup (\tetrois \teun^i \tedeux^j)_{i,j \in \mathbb{N}}$ is a basis of $A_{\alpha,\beta}$ viewed as a $\mathbb{C}$-vector space. So we can write $a$ as follows:
$$a= \sum_{(i,j) \in I} \alpha_{i,j} \teun^i \tedeux^j + \sum_{(i,j) \in I'} \beta_{i,j} \tetrois \teun^i \tedeux^j,$$
where $I,I'$ are two finite subsets of $\mathbb{N}^2$, $\alpha_{i,j} \in \mathbb{C}$ for all $(i,j) \in I$ 
and $\beta_{i,j} \in \mathbb{C}$ for all $(i,j) \in I'$. Thus, since $\teun a=1$  and $\teun \tetrois = q^{-2} \tetrois \teun$, we get : 
$$ \sum_{(i,j) \in I} \alpha_{i,j} \teun^{i+1} \tedeux^j + \sum_{(i,j) \in I'} \beta_{i,j}q^{-2} \tetrois \teun^{i+1} \tedeux^j=1.$$
Since $\mathcal{F}$ is linearly independent, we obtain that $\alpha_{i,j} =0$ for all $(i,j) \in I$ 
and $\beta_{i,j} =0$ for all $(i,j) \in I'$, so that $a=0$. This is a
contradiction and so we have just proved that $\teun$ is not
invertible in $A_{\alpha,\beta}$.

Assume now that $\tetrois $ is invertible in $A_{\alpha,\beta}$. Then there exists $a \in A_{\alpha,\beta}$ such that $\tetrois a =1$. As in the previous case, we can write $a$ as follows:
$$a= \sum_{(i,j) \in I} \alpha_{i,j} \teun^i \tedeux^j + \sum_{(i,j) \in I'} \beta_{i,j} \tetrois \teun^i \tedeux^j,$$
where $I,I'$ are two finite subsets of $\mathbb{N}^2$, $\alpha_{i,j} \in \mathbb{C}$ for all $(i,j) \in I$ 
and $\beta_{i,j} \in \mathbb{C}$ for all $(i,j) \in I'$. Thus, since $\tetrois a=1$ and $\teun \tetrois = q^{-2} \tetrois \teun$, we get:
$$ \sum_{(i,j) \in I} \alpha_{i,j} \tetrois \teun^{i} \tedeux^j + \sum_{(i,j) \in I'} \beta_{i,j}q^{4i}  \teun^{i} \tetrois^2 \tedeux^j=1.$$
Now recall (see \ref{soussectionbase}) that $\tetrois^2+c_1 \tetrois \teun \tedeux +\alpha c_2 \teun + \beta c_3 =0$, where $c_1= q^4-1$, $c_2=q^2(q^2+1)$ and $c_3=\frac{q^6}{1-q^2}$. Hence we deduce from the previous equality that 
\begin{eqnarray*}
 \sum_{(i,j) \in I} \alpha_{i,j} \tetrois \teun^i \tedeux^j  & - &   c_1\sum_{(i,j) \in I'} q^{2i} \beta_{i,j} \tetrois \teun^{i+1} \tedeux^{j+1} - \: c_2 \sum_{(i,j) \in I'} q^{4i} \alpha \beta_{i,j} \teun^{i+1}  \tedeux^j  \\
 - \: c_3 \sum_{(i,j) \in I'} q^{4i} \beta \beta_{i,j} \teun^{i}  \tedeux^j & = & 1  \\
\end{eqnarray*}
If every $\beta_{i,j}$ ($(i,j) \in I'$) is equal to zero, then $\sum_{(i,j) \in I} \alpha_{i,j} \tetrois \teun^i \tedeux^j =1$. Since the family $\mathcal{F}$ is linearly independent, this implies that $\alpha_{i,j} =0$ for all $(i,j) \in I$. Hence $a =0$. This is a contradiction and so  there exists 
$(u,v) \in I'$ such that $\beta_{u,v} \neq 0$. Choose such a pair $(u,v)$ maximal with respect to the lexicographic order. Since $\mathcal{F}$ is linearly independent, identifying the coefficients of $ \teun^{u+1} \tedeux^v$ in the previous equality  leads to $$q^{4u}\alpha \beta_{u,v} = 0.$$
This is a contradiction since $\beta_{u,v}  \neq 0$ and  $\alpha
\neq 0$.\\$ $ 

To sum up, $\teun$ and $\tetrois$ are not invertible in
$A_{\alpha,\beta}$.\\$ $

\noindent $\bullet$ Step 2: the case where $\beta \neq 0$. \\

Let  $u \in U(A_{\alpha,\beta})$. It follows from the discussion
preceding Proposition \ref{unite} that there exist $m,n \in \mathbb{Z}$ and
$\lambda \in \mathbb{C}^*$ such that $u= \lambda \tetrois^m
\teun^n$. We need to prove that $m=n=0$.

First, observe that, since $u^{-1}=\lambda^{-1}q^{-2mn}\tetrois^{-m}
\teun^{-n}$ also belongs to $U(A_{\alpha,\beta})$, we can assume that
$m \geq 0$.

Now, since  $\teun$ and $\tetrois$ are not invertible in
$A_{\alpha,\beta}$, it is easy to show that we must have 
\\either $m=n=0$ (and $u \in \mathbb{C}^*$ in this case), or $m>0$ and $n<0$.

Let us study this last case, that is, assume that $m>0$ and $n<0$. We set $k:=-n >0$. 
Then $u=\lambda \tetrois^m \teun^{-k}=\lambda q^{\bullet} \teun^{-k}
\tetrois^m$, where $\bullet \in \mathbb{Z}$. Thus $\lambda q^{\bullet} \tetrois^m= \teun^k u$. Now $u \in 
A_{\alpha,\beta}$, so that one can write $u$ uniquely as follows:
$$u= \sum_{(i,j) \in I} \alpha_{i,j} \teun^i \tedeux^j + \sum_{(i,j) \in I'} \beta_{i,j} \tetrois \teun^i \tedeux^j,$$
where $I,I'$ are two finite subsets of $\mathbb{N}^2$, $\alpha_{i,j}
\in \mathbb{C}$ for all $(i,j) \in I$ and $\beta_{i,j}
\in \mathbb{C}$ for all $(i,j) \in I'$. Hence we get: 
\begin{eqnarray*}
\lambda q^{\bullet} \tetrois^m= \sum_{(i,j) \in I} \alpha_{i,j} \teun^{i+k}
\tedeux^j + \sum_{(i,j) \in I'}q^{-2k} \beta_{i,j} \tetrois
\teun^{i+k} \tedeux^j.
\end{eqnarray*}
Since $\beta \neq 0$, comparing this expression of $ \tetrois^m$ 
to the expression of $ \tetrois^m$ obtained in Lemma 
\ref{e3base} leads to a contradiction. Hence the only possibility is
$m=n=0$, so that $u=\lambda \in \mathbb{C}^*$ and
$U(A_{\alpha,\beta})=\mathbb{C}^*$, as desired.\\$ $

\noindent $\bullet$ Step 3: the case where $\beta = 0$. \\

Let  $u \in U(A_{\alpha,0})$. It follows from the discussion
preceding Proposition \ref{unite} that there exist $m,n \in \mathbb{Z}$ and
$\lambda \in \mathbb{C}^*$ such that $u= \lambda \tetrois^m
\teun^n$. We need to prove that $m=-2n$.

First, as in the previous step, one can assume that $m \geq 0$.

Now, since  $\teun$ and $\tetrois$ are not invertible in
$A_{\alpha,0}$, it is easy to show that we must have 
\\either $m=n=0$ (and $u \in \mathbb{C}^*$ in this case), or $m>0$ and $n<0$.  

Assume that $m>0$ and $n<0$. As in the second step, we set $k:=-n
>0$. Keeping the notations of the previous step, one can also prove
that 
\begin{eqnarray}
\label{tototo}
\lambda q^{\bullet} \tetrois^m= \sum_{(i,j) \in I} \alpha_{i,j} \teun^{i+k}
\tedeux^j + \sum_{(i,j) \in I'}q^{-2k} \beta_{i,j} \tetrois
\teun^{i+k} \tedeux^j.
\end{eqnarray}
Then, comparing (\ref{tototo}) to the expression of $ \tetrois^m$ obtained in Lemma 
\ref{e3base}, we get that $2k \leq m$. Thus we can rewrite $u$ as follows:
\begin{eqnarray}
\label{uuu}
u = \lambda q^{*} \tetrois^{m-2k} \left(\teun^{-1} \tetrois^2 
\right)^k,
\end{eqnarray}
where $* \in \mathbb{Z}$. Observe further that, since $\beta=0$, we have $\tetrois^2+c_1
\tetrois \teun \tedeux +\alpha c_2 \teun  =0$. Thus 
$$\teun^{-1} \tetrois^2 = -q^2 c_1 \tetrois \tedeux -\alpha c_2 \in A_{\alpha,0}.$$
So, recalling that $2k \leq m$, (\ref{uuu}) expresses $u$ as a product of elements of
$A_{\alpha,0}$. Since $u$ is invertible, this implies that all
factors of this expression should  also be invertible in
$A_{\alpha,0}$. However we have already proved that
$\tetrois$ is not invertible in $A_{\alpha,0}$. So this implies that
$m=2k$, so that $u = \lambda q^{*}  \left(\teun^{-1} \tetrois^2 
\right)^k$. 
\\$ $

Hence we have already proved that 
$U(A_{\alpha,0}) \subseteq \{\lambda (\teun^{-1} \tetrois^2)^i \mid
\lambda \in \mathbb{C}^* \mbox{, } i \in \mathbb{Z} \}$. 
\\$ $

For the converse, it is sufficient to observe that $\teun^{-1}
\tetrois^2$ is invertible in $A_{\alpha,0}$ with inverse:
$$(\teun^{-1}\tetrois^2 )^{-1}= \frac{1}{\alpha^2 c_2^2} \left[ -\alpha
  c_2 +c_1^2 \teun \tedeux^2 +c_1 \tetrois \tedeux   \right]. \mbox{ \fin} $$
$ $

If $A$ is an algebra, the trivial invertible elements of $A$ are the non-zero elements of the ground field.

\begin{prop}
Let $ \alpha, \beta \in \mathbb{C}^*$.
\\Then $A_{\alpha,\beta}$ is a simple algebra with no non-trivial invertible elements and the center of $A_{\alpha,\beta}$ is reduced to $\mathbb{C}$.
\end{prop}
\preuve It just remains to show that the center of $A_{\alpha,\beta}$ is reduced to $\mathbb{C}$. 
Let $x$ a non-zero element of the center of $A_{\alpha,\beta}$. Then $x A_{\alpha,\beta}$ is a two-sided non-zero ideal 
of $A_{\alpha,\beta}$. Since this algebra is simple, we get
$xA_{\alpha,\beta}=A_{\alpha,\beta}$. This implies that $x$ is
invertible in $A_{\alpha,\beta}$. Then Proposition \ref{unite} leads to $x \in \mathbb{C}$, as
desired. \fin
\\$ $

One can also show, using similar arguments, that the center of the algebra $A_{\alpha,\beta}$ with $\alpha
\beta =0$ is also reduced to $\mathbb{C}$.

\begin{prop}
Let $\alpha, \beta \in \mathbb{C}^*$.
\\Then $A_{0,\beta}$ and  $A_{\alpha,0}$ are simple algebras whose centers are reduced to $\mathbb{C}$.
\end{prop}

Note that, in the case of the so-called
Weyl-Hayashi algebras $A_{0,\beta}$ ($\beta \in \mathbb{C}^*$), this result was first proved by Kirkman and Small (see \cite[Proposition 1.5]{kirkman}).

\subsection{Link with generalized Weyl algebras.}
\label{GWA}
$ $

It is interesting to note that the algebras $A_{0,\beta}$ ($\beta \neq
0$) and $A_{\alpha,0}$ ($\alpha \neq 0$) are examples of a large class of algebras 
called generalized Weyl algebras (GWA for short) that have been introduced by Bavula (see \cite{bavula}). These algebras have been extensively studied recently. We refer to \cite{bavula,bavjor,lionel} for more details on these algebras and their representation theory. Recall that a GWA 
over a Noetherian $\mathbb{C}$-algebra $R$ which is a domain is defined as follows. Let $a$ be a non-zero central element of $R$ and $\sigma $ be an automorphism of $R$. Then the GWA $R(\sigma,a)$ is the $\mathbb{C}$-algebra generated over $R$ by two generators $x$ and $y$ with the following relations:
$$xr=\sigma(r) x \mbox{ for all } r \in R,$$
$$yr=\sigma^{-1}(r) y \mbox{ for all } r \in R,$$
$$xy= \sigma(a),$$
$$yx = a.$$
The algebra $R(\sigma,a)$ is a Noetherian domain (see \cite{bavula}).

In fact, the algebras $A_{0,\beta}$ ($\beta \neq 0$) and  $A_{\alpha,0}$ ($\alpha \neq 0$) are examples of GWA over a Laurent polynomial ring in one indeterminate since it can be shown that:

\begin{prop}
\label{gwa1}
Let $\beta \in \mathbb{C}^*$. Denote by $\sigma$ the automorphism of $\mathbb{C}[h^{ \pm 1}]$ defined 
by $\sigma(h) = q^2 h$ and set $a= \frac{h}{1-q^4}+ \frac{\beta q^6}{(q^2-1)(q^4-1)}h^{-1}$.
\\Then there exists a unique isomorphism $\Theta : A_{0,\beta} \rightarrow \mathbb{C}[h^{ \pm 1}]( \sigma, a)$ such that 
$\Theta (\teun ) = y$, $\Theta (\tedeux ) = x$ and $\Theta (\tetrois ) = h$. 
\end{prop}
and 
\begin{prop}
\label{gwa2}
Let $\alpha \in \mathbb{C}^*$. Denote by $\sigma$ the automorphism of $\mathbb{C}[h^{ \pm 1}]$ defined 
by $\sigma(h) = q^2 h$ and set $a= \frac{h^{-1}}{1-q^4}- \frac{\alpha}{q^2-1}$.
\\Then there exists a unique isomorphism $\Theta : A_{\alpha,0} \rightarrow \mathbb{C}[h^{ \pm 1}]( \sigma, a)$ such that 
$\Theta (\teun ) = x^2 h$, $\Theta (\tedeux ) = y$ and $\Theta (\tetrois ) = x$. 
\end{prop}

Let $\alpha,\beta \in \mathbb{C}^*$. In contrast, the algebra $A_{\alpha,\beta}$ is not isomorphic to a GWA over 
$\mathbb{C}[h]$ nor $\mathbb{C}[h^{\pm 1}]$. In fact, every GWA over
$\mathbb{C}[h^{\pm 1}]$ has non-trivial invertible elements ($h$ is
invertible in this algebra) and so Proposition \ref{unite} shows that
$A_{\alpha,\beta}$ cannot be isomorphic to a GWA over
$\mathbb{C}[h^{\pm 1}]$. Now, assume that $A_{\alpha,\beta}$ is
isomorphic to a GWA over $\mathbb{C}[h]$. Since the field of fractions
of $A_{\alpha,\beta}$ is isomorphic to the field of fractions of a
quantum plane, then it follows from \cite[Proposition 2.1.1]{lionel}
that $A_{\alpha,\beta}$ must be isomorphic to a quantum GWA over
$\mathbb{C}[h]$ (in the sense of \cite{lionel}). But, it was observed
in \cite[3.1]{lionel} that such quantum GWA are never simple. This is a contradiction and so we have just proved the following statement.

\begin{prop}
Let $(\alpha,\beta) \in \mathbb{C}^2$ with $\alpha,\beta \neq 0$.
\\Then $A_{\alpha,\beta}$ is not isomorphic to a GWA over 
$\mathbb{C}[h]$ nor $\mathbb{C}[h^{\pm 1}]$.
\end{prop}

\section{Automorphism group of $U^+$ and orbits of the action of this group on $\prim(U^+)$.}
$ $

The aim of this paragraph is to compute the automorphism group $\aut (U^+)$ of the algebra $U^+$ and to describe the orbits for the action of $\aut (U^+)$ on $\prim(U^+)$. In \cite{andrusdumas}, Andruskiewitsch and Dumas have investigated the group $\aut (U^+)$ and conjectured that this group is isomorphic to the torus $\hc=(\mathbb{C}^*)^2$ (see \cite[Problem 1]{andrusdumas}). To give support to their conjecture, they have proved the following intermediate result (see \cite[Proposition 3.3]{andrusdumas}). 

\begin{prop}[Andruskiewitsch-Dumas]
\label{dumas}
The sub-group $\aut_z(U^+)$ of $\aut (U^+)$ of those automorphisms that fix the ideal generated by $z$ is isomorphic to the torus $(\mathbb{C}^*)^2$. 
\end{prop}

 Our aim in this section is to give a positive answer to the conjecture of Andruskiewitsch and Dumas. In fact, in view of Proposition \ref{dumas}, it is sufficient to prove that every automorphism of $U^+$ fixes the ideal generated by $z$. In order to do this, we proceed as follows. Using the results obtained in the previous paragraph, we establish that the set $\prim^z(U^+)$ of those primitive ideals of $U^+$ that contain $z$ is left invariant by every automorphism of $U^+$. Next, since $U^+$ is a Jacobson ring, every prime ideal is the intersection of those primitive ideals that contain this prime. In particular, the ideal generated by $z$ is the intersection of the set $\prim^z(U^+)$ of those primitive ideals that contain $z$. Since we have previously proved that the set $\prim^z(U^+)$ is invariant under any automorphism, we conclude that the ideal generated by $z$ is also invariant under any automorphism. Then Proposition \ref{dumas} of Andruskiewitsch and Dumas allows us to prove that the group $\aut(U^+)$ is isomorphic to the torus $(\mathbb{C}^*)^2$. As a corollary, we obtain that the action of $\aut (U^+)$ on $\prim(U^+)$ has exactly $8$ orbits that we describe explicitly.

\subsection{Image of the maximal ideals $\langle z ,z'-\beta \rangle$ ($\beta \in \mathbb{C}^*$) by an automorphism of $U^+$.}
$ $

Throughout this section, $\beta$ denotes a non-zero complex number and $\sigma $ an automorphism of $U^+$. Further, we keep the conventions and notations of \ref{soussectionbase}. 
Since $\langle z ,z'-\beta \rangle$ is a maximal ideal of height 2 of $U^+$, its image under $\sigma$ is also a maximal ideal of height 2 of $U^+$. So we deduce from Proposition \ref{max} that there exists $(\alpha',\beta') \in \mathbb{C}^2 \setminus \{(0,0)\}$ such that $\sigma \left( \langle z ,z'-\beta \rangle \right) = \langle z-\alpha' ,z'-\beta' \rangle$. Further, since $\sigma$ induces an isomorphism between $U^+ / \langle z ,z'-\beta \rangle$ and 
$U^+ / \langle z-\alpha' ,z'-\beta' \rangle$, the algebras
$A_{0,\beta}$ and $A_{\alpha',\beta'}$ must be isomorphic. In
particular, their groups of invertible elements must be
isomorphic. Hence, it follows from Proposition \ref{unite} that we
must have $\alpha'=0$ or $\beta'=0$. 

Next, let $\alpha \in \mathbb{C}^*$. Even though $A_{0,\beta}$ and $A_{\alpha,0}$ are isomorphic to GWA over a Laurent
polynomial ring in one variable (see Proposition \ref{gwa1} and
\ref{gwa2}), they are not isomorphic algebras. Indeed, it follows from \cite[Theorem 5.2]{bavjor} that the GWA that are
isomorphic to $A_{0,\beta}$ and $A_{\alpha,0}$ (see Proposition \ref{gwa1} and
\ref{gwa2}) are not isomorphic.

Hence, we must have $\alpha'=0$ and so we have just proved the following result.


\begin{prop}
\label{separation}
 Let $\beta \in \mathbb{C}^*$ and $\sigma $ an automorphism of $U^+$. 
\\There exists $\beta' \in \mathbb{C}^*$ such that $\sigma \left( \langle z ,z'-\beta \rangle \right) = \langle z,z'-\beta' \rangle$.
\end{prop}

\subsection{Automorphism group of $U^+$.}
$ $

We denote by $\aut(U^+)$ the automorphism group of the algebra $U^+$ and by $\aut_{z}(U^+)$ the sub-group of $\aut(U^+)$ of those automorphisms of $U^+$ that fix $\langle z \rangle$. Recall that Andruskiewitsch and Dumas have shown (see \cite[Proposition 3.3]{andrusdumas}) that $\aut_{z}(U^+)=\{\psi_{\alpha,\beta} \mid \alpha, \beta \in \mathbb{C}^* \}$, where for all $\alpha, \beta \in \mathbb{C}^*$, $\psi_{\alpha,\beta}$ denotes the automorphism of $U^+$ defined by:
$$ \psi_{\alpha,\beta}(e_1)=\alpha e_1 \mbox{ and }\psi_{\alpha,\beta}(e_2)= \beta e_2.$$
Hence, $\aut_{z}(U^+)$ is isomorphic to the torus $(\mathbb{C}^*)^2$. In order to prove that $\aut(U^+)$ itself is isomorphic to $(\mathbb{C}^*)^2$, it just remains to prove that $\aut_{z}(U^+)= \aut(U^+)$. This is what we do now.


\begin{prop}
\label{fixe}
  $\aut_{z}(U^+)= \aut(U^+)$.
\end{prop}
\preuve The inclusion $\aut_{z}(U^+) \subseteq  \aut(U^+)$ is obvious. Let now $\sigma$ be an automorphism of $U^+$. In order to prove that the ideal generated by $z$ is left invariant by $\sigma$, we first show that the set $\prim^z(U^+)$ of those primitive ideals of $U^+$ that contain $z$ is invariant under $\sigma$.

It follows from Proposition \ref{prim} that $\prim^z(U^+)$ is actually the disjoint union of the following $\hc$-strata of $\prim(U^+)$: $ \{ \langle z ,z'-\beta \rangle \mid \beta \in \mathbb{C}^* \} $, 
$ \{ \langle e_3 \rangle \} $, 
$\{ \langle \ov{e_3} \rangle \}$, 
$\{ \langle e_1, e_2-\beta \rangle  \mid \beta \in \mathbb{C}^* \} $, 
$ \{ \langle e_1 - \alpha , e_2 \rangle  \mid \alpha \in \mathbb{C}^* \} $, 
and $\{ \langle e_1, e_2 \rangle \}$. 
We set $\mathcal{P}_z:=\{ \langle z ,z'-\beta \rangle \mid \beta \in \mathbb{C}^* \}$, 
$\mathcal{P}_0 := \{ \langle e_3 \rangle \} \:  \bigsqcup  \:  \{ \langle \ov{e_3} \rangle \}$ 
and $\prim^4(U^+) =  \{ \langle e_1, e_2-\beta \rangle  \mid \beta \in \mathbb{C}^* \} \: \bigsqcup  \:   \{ \langle e_1 - \alpha , e_2 \rangle  \mid \alpha \in \mathbb{C}^* \} \: \bigsqcup  \:   \{ \langle e_1, e_2 \rangle \}$, so that $\prim^z(U^+)= \mathcal{P}_z \: \bigsqcup \: \mathcal{P}_0 \: \bigsqcup \: \prim^4(U^+)$. Now, to prove that $\prim^z(U^+)$ is invariant under the automorphism $\sigma$, it is sufficient to show that  the three subsets $\mathcal{P}_z$, $\mathcal{P}_0$ and $\prim^4(U^+)$ are left invariant by $\sigma$.

First, it follows from Proposition \ref{separation} that the image by $\sigma$ of an element of 
$\mathcal{P}_z$ is still an element of $\mathcal{P}_z$. Hence, $\mathcal{P}_z$ is invariant under $\sigma$. Next, it follows from Proposition \ref{max} that $\mathcal{P}_0$ is exactly the set of those primitive ideals of $U^+$ that are not maximal. Thus $\mathcal{P}_0$ is also invariant under $\sigma$. Finally, it follows from Proposition \ref{prim} that $\prim^4(U^+)$ is actually the set of all height $4$ prime ideals of $U^+$. Hence this set is also invariant under $\sigma$. To resume, we have just proved that the three sets $\mathcal{P}_z$, $\mathcal{P}_0$ and $\prim^4(U^+)$ are invariant under $\sigma$. Since $\prim^z(U^+)= \mathcal{P}_z \: \bigsqcup \: \mathcal{P}_0 \: \bigsqcup \: \prim^4(U^+)$, we obtain that $\prim^z(U^+)$ is invariant under $\sigma$ too, that is, $\sigma \left( \prim^z(U^+) \right) = \prim^z(U^+)$. 

 Next, it follows from \cite[Proposition II.7.12]{bg} that $U^+$ is a Jacobson ring. In particular, this implies that the prime $\langle z \rangle$ is actually equal to the intersection of the elements of $\prim^z(U^+)$, that is, the intersection of those primitive ideals of $U^+$ that contain $z$. In other words, we have
$$ \langle z \rangle = \bigcap_{P \in \prim^z(U^+) } P.$$ 
Applying the automorphism $\sigma$ to this equality yields
$$ \sigma \left(  \langle z \rangle \right)  = \bigcap_{P \in \prim^z(U^+) } \sigma(P).$$
Since $\prim^z(U^+)$ is left invariant by $\sigma$, this implies that  
 $$ \sigma \left(  \langle z \rangle \right)  = \bigcap_{P \in \prim^z(U^+) } P = \langle z \rangle,$$
as desired. \fin \\$ $


Since Andruskiewitsch and Dumas have shown that  $\aut_{z}(U^+)$ is isomorphic to the torus $(\mathbb{C}^*)^2$ 
(see \cite[Proposition 3.3]{andrusdumas}), we deduce from Proposition \ref{fixe} the following result that give a positive answer to a conjecture of Andruskiewitsch and Dumas (see \cite[Problem 1]{andrusdumas}).

\begin{theo}
\label{thetheoauto}
$\aut(U^+)=\{\psi_{\alpha,\beta} \mid \alpha, \beta \in \mathbb{C}^* \}$, where, for all $\alpha, \beta \in \mathbb{C}^*$, $\psi_{\alpha,\beta}$ denotes the automorphism of $U^+$ defined by $ \psi_{\alpha,\beta}(e_1)=\alpha e_1 \mbox{ and }\psi_{\alpha,\beta}(e_2)= \beta e_2$.
\\Hence, $\aut(U^+)$ is isomorphic to the torus $(\mathbb{C}^*)^2$.
\end{theo} 

\subsection{Orbits for the action of $\aut(U^+)$ on $\prim(U^+)$. }
$ $

We easily deduce from \cite[Theorem II.8.14]{bg} and Theorem \ref{thetheoauto} the following statement.

\begin{cor}
The $\aut(U^+)$-orbits within $\prim(U^+)$ coincide with the $\hc$-strata of $\prim(U^+)$. 
\\In other words, the action of $\aut(U^+)$ on $\prim(U^+)$ has exactly $8$ orbits:
\begin{enumerate}
\item $\prim_{\langle 0 \rangle } (U^+)  =  \{ \langle z-\alpha ,z'-\beta \rangle \mid \alpha,\beta \in \mathbb{C}^* \} $,
\item $\prim_{\langle z' \rangle } (U^+)  =   \{ \langle z-\alpha ,z' \rangle \mid \alpha \in \mathbb{C}^* \}$ ,
\item $\prim_{\langle z \rangle } (U^+)  =  \{ \langle z ,z'-\beta \rangle \mid \beta \in \mathbb{C}^* \} $,
\item $\prim_{\langle e_3 \rangle } (U^+)   =  \{ \langle e_3 \rangle \}$,
\item $\prim_{\langle \ov{e_3} \rangle } (U^+)  =  \{ \langle \ov{e_3} \rangle \}$,
\item $\prim_{\langle e_1 \rangle } (U^+)   =   \{ \langle e_1,e_2-\beta \rangle \mid \beta \in \mathbb{C}^* \}$,
\item $\prim_{\langle e_2 \rangle } (U^+)  =  \{ \langle e_1-\alpha , e_2 \rangle \mid \alpha \in \mathbb{C}^* \}$, 
\item $\prim_{\langle e_1,e_2 \rangle } (U^+)  =  \{ \langle e_1,e_2 \rangle \}$.
\end{enumerate}
\end{cor}

%
%

\noindent {\it Current contact details (from October 2005):}
\\St\'ephane Launois, 
\\University of Edinburgh, School of Mathematics 
\\James Clerk Maxwell Building, Kings Buildings, Mayfield Road 
\\Edinburgh EH9 3JZ
\\ United Kingdom.
\\email: stephane.launois@ed.ac.uk

\end{document}